\documentclass[10pt]{article}
\usepackage{amsfonts}
\usepackage{mathrsfs}
\usepackage{amsmath,amsthm,amsfonts,amssymb,epsfig}

\newtheorem{thm}{Theorem}[section]

\newtheorem{lem}[thm]{Lemma}

\theoremstyle{definition}
\newtheorem{defn}[thm]{Definition}
\theoremstyle{remark}

\numberwithin{equation}{section}

\newcommand\be{\begin{equation}}
\newcommand\ee{\end{equation}}
\newcommand\bes{\begin{eqnarray}}
\newcommand\ees{\end{eqnarray}}
\newcommand\bess{\begin{eqnarray*}}
\newcommand\eess{\end{eqnarray*}}

\begin{document}
\title{Null controllability and inverse source problem for stochastic Grushin equation with boundary degeneracy and singularity}
\author{ \ Lin Yan$^a$, \ Bin Wu$^a$\footnote{Corresponding author. email: binwu@nuist.edu.cn},  \ Shiping Lu$^a$,\ Yuchan Wang$^a$  \\
$ ^a$School of Mathematics and Statistics\\ Nanjing University of
Information Science
and Technology, \\
Nanjing 210044, P. R. China
}

\maketitle

\begin{abstract}
In this paper, we consider a null controllability and an inverse source problem for stochastic Grushin equation with boundary degeneracy and singularity. We construct two special weight functions to establish two Carleman estimates for the whole stochastic Grushin operator with singular potential by a weighted identity method. One is for the backward stochastic Grushin equation with singular weight function. We then apply it to prove the null controllability for stochastic Grushin equation for any $T$ and any degeneracy $\gamma>0$, when our control domain touches the degeneracy line $\{x=0\}$. In order to study the inverse source problem of determining two kinds of sources simultaneously, we prove the other Carleman estimate, which is for the forward stochastic Grushin equation with regular weight function. Based on this Carleman estimate, we obtain the uniqueness of the inverse source problem.
\vskip 0.3cm


\vskip 0.3cm


\end{abstract}
\section{Introduction}
\setcounter{equation}{0}

Let $(\Omega,\mathcal F, \{\mathcal F_t\}_{t\geq0},\mathbb P)$ be a complete filtered probability space, on which a one-dimensional standard Brownian motion $\{B(t)\}_{t\geq0}$ is defined.
Let $I=I_x\times I_y$ with $I_x=(0,1)$,
$I_y=(0,1)$, $Q_T= I\times(0,T)$, $\Sigma_T=\partial I\times(0,T)$. Then we consider the following stochastic
Grushin equation with singular potential:
\bes\label{1.1}\left\{\begin{array}{ll}
{\rm d}u-u_{xx}{\rm d}t-x^{2\gamma}u_{yy}{\rm d}t-\frac{\sigma}{x^2}u{\rm d}t=f{\rm d}t+F{\rm d}B(t),&(x,y,t)\in Q_T,\\
u(x,y,t)=0,&(x,y,t)\in \Sigma_T,\\
u(x,y,0)=u_0(x,y),&(x,y)\in I,
\end{array}\right.
\ees
where $\sigma$ and $\gamma$ are  two constants. Obviously, the system (\ref{1.1}) is not only degenerate, but also singular on boundary $\{x=0\}\times I_y$. Further, the degeneracy is weak if $0<\gamma<\frac{1}{2}$ and strong if $\gamma\geq  \frac{1}{2}$.

This paper focus on the Carleman estimates for stochastic
Grushin equation with singular potential and then apply them to
study the following null controllability and  inverse  source problem.

Here and henceforth, for any $a\in (0,1)$ we set $\nabla=(\partial_x,\partial_y)$ and
\begin{align*}\begin{array}{ll} \omega=(0,a)\times I_y,&\omega_T=\omega\times (0,T),\\
               \Gamma=\{x=0\}\times I_y,&\Gamma_T=\Gamma\times (0,T),
               \end{array}
\end{align*}
where $\omega$ is the control domain for null controllability, $\Gamma$ is the observation boundary for inverse source problem. It is noted that our control domain touches the degeneracy line $\{x=0\}$ as [\ref{Beauchard2015}], where the null controllability for the deterministic Grushin equation without singularity, i.e. $\sigma=0$,  is obtained for any $T$ and any $\gamma>0$.
\vspace{2mm}

{\noindent\bf Null Controllability.}\
For any $u_0\in L^2(\Omega,\mathcal F_0, \mathbb P; L^2(I))$, find a pair  $(g, G)$ such that the solution $u$ of the following forward stochastic
Grushin equation with singular potential:
\bes\label{1.2}\left\{\begin{array}{ll}
{\rm d}u-u_{xx}{\rm d}t-x^{2\gamma}u_{yy}{\rm d}t-\frac{\sigma}{x^2}u{\rm d}t=(\alpha u+g{\bf 1}_\omega){\rm d}t+(\beta u+G){\rm d}B(t),&(x,y,t)\in Q_T,\\
u(x,y,t)=0,&(x,y,t)\in \Sigma_T,\\
u(x,y,0)=u_0(x,y),&(x,y)\in I,
\end{array}\right.
\ees
satisfies
\begin{align*}u(x,y,T)=0,\quad  (x,y)\in I,\  \mathbb P-a.s.,
\end{align*}
where ${\bf 1}_\omega$ is the characteristic function of the set $\omega$.

\vspace{2mm}

{\noindent\bf Inverse source problem.}\
Determine two kinds of sources  $h(x,t)$ and $H(t)$ simultaneously in the following forward stochastic
Grushin equation with singular potential:
\begin{equation}\label{1.3}
	\left\{
	\begin{aligned}
		&
		\begin{aligned}
		{\rm d}u-u_{xx}{\rm d}t-x^{2\gamma}u_{yy}{\rm d}t-\frac{\sigma}{x^2}u{\rm d}t =&h(x,t)R_1(x,y,t){\rm d}t \\
		& +H(t)R_2(x,y,t){\rm d}B(t), \quad (x,y,t)\in Q_T,
		\end{aligned} \\
		& u(x,y,t)=0, \quad\quad\quad\qquad\qquad\qquad\qquad\qquad\qquad\qquad\qquad\quad \ \ (x,y,t)\in\Sigma_T,\\
        & u(x,y,0)=0, \quad\quad\quad\qquad\qquad\qquad\qquad\qquad\qquad\qquad\qquad\quad \ \ (x,y)\in I,\\
	\end{aligned}
	\right.
\end{equation}
by the boundary observation $u_y|_{\Sigma_T}$, $u_x|_{\Gamma_T}$ and final time observation $u|_{t=T}$ in $I$.

\vspace{2mm}

When no singular term was involved, the null controllability of deterministic Grushin equation with $I=(-1,1)\times (0,1)$ was
studied in [\ref{Anh2013},\ref{Beauchard2014}].  The null controllability for Grushin-type equations was obtained for any time $T>0$ and for any degeneracy $\gamma>0$, with a control that acts on one strip, touching the degeneracy line $\{x = 0\}$ in [\ref{Beauchard2015}].  When restricting the domain to one side only of the singular set, i.e. $I=(0,1)\times (0,1)$, [\ref{Cannarsa2014}] proved that there exists $T^*$ such that for every $T>T^*$ the Grushin-type equation is null
controllable for $\gamma=1$, $\sigma<\frac{1}{4}$. Next, [\ref{Anh2016}] showed a similar null controllability in large time $T$ when the degeneracy of the diffusion coefficient
and singularity of the potential occur at the interior of the domain. The key ingredient in these papers is applying a Fourier decomposition to reduce the problem to the validity of a uniform observability inequality with respect to the Fourier frequency. As for the inverse source problem for deterministic Grushin equation, [\ref{Beauchard2014IP}] proved a Lipschitz stability result of determining a source function $h$ depending on $x$ and $y$,  by the observation data $\partial_tu|_{\omega\times(T_1,T_2)}$ with a suitable subdomain $\omega$.

It is well known that Carleman estimate is the key tool to study null controllability and inverse problems, which is a class of  weighted energy estimates in connection with deterministic/stochastic differential operators. As its applications to deterministic differential equations, we refer to [\ref{Cannarsa2010},\ref{Jiang2017},\ref{Klibanov2004},\ref{Klibanov2013},\ref{Wu2017},\ref{Yamamoto2009}] for inverse problems, [\ref{Buhgeim1981},\ref{Rousseau2012},\ref{Saut1987},\ref{Wu2019JIIP}] for unique continuation problems, [\ref{Imanuvilov2003},\ref{Gao2016},\ref{lu2013},\ref{Fragnelli2016}] for control theory. For Carleman estimates related to deterministic Grushin equation, we refer to [\ref{Anh2013},\ref{Beauchard2014},\ref{Morancey2013},\ref{Koenig}]. In recent years, many efforts have been devoted to studying the Carleman estimate for stochastic partial differential equations, for example   [\ref{Barbu2003},\ref{Liu2014},\ref{Zhang2009},\ref{Yan2018JMAA}] for stochastic heat equation, [\ref{Zhang2008}] for stochastic wave equation, [\ref{Gao2014}] for stochastic KdV equation, [\ref{Gao2015}] for stochastic Kuramoto-Sivashinsky equation, [\ref{2013}] for stochastic Schr\"{o}dinger equation, and so on. To the best of our knowledge, there is only one paper about Carleman estimates for one dimensional stochastic
degenerate operator ${\rm d}u-x^{2\gamma} u_{xx}{\rm d}t$ [\ref{LiuSIAM2019}], which is very different from the degenerate Grushin operator ${\rm d}u-u_{xx}{\rm d}t-x^{2\gamma}u_{yy}{\rm d}t$. In these works, Carleman estimates were mainly applied to deal with stochastic control problems.
Since the solution of a stochastic differential equation is not differentiable with respect to time variable, which leads to that some traditional methods for deterministic inverse problems  cannot be applied to the corresponding ones in the stochastic case. Therefore, [\ref{Lu2015CPAM}] proposed a regular weight function in Carleman estimates to study an stochastic inverse problem related to the stochastic hyperbolic equation. We also refer to [\ref{2012},\ref{Yuan2015}]  for stochastic inverse problems.

 Although there are numerous results for Carleman estimates for deterministic Grushin equation,  little has been known for Carleman estimates related to the stochastic Grushin equation. In this paper, we  first construct a special weight function $\psi$ to obtain a Carleman estimate for backward stochastic Grushin operator with singular potential and then apply this Carleman estimate to prove the null controllability for system \eqref{1.2}.  We do not apply the method based on Fourier decomposition as [\ref{Anh2016},\ref{Cannarsa2014}].  A weakness of Fourier decomposition is that in proving the observation inequality the authors have to deal with the eigenvalues in Fourier decomposition $\mu_{n}\rightarrow +\infty$ as $n\rightarrow\infty$, which is the reason that the condition $T>T^*$ is introduced in [\ref{Cannarsa2014}]. In order to obtain the null controllability result for any time $T$ and any degeneracy $\gamma$, we consider the Grushin operator with singular potential, i.e $u_{xx}+x^{2\gamma}u_{yy}+\frac{\sigma}{x^{2}}$, as a whole to establish our Carleman estimate, not as [\ref{Cannarsa2014}] only for its Fourier components with respect to $u$, i.e. $(u_n)_{xx}-\big[(n\pi)^2x^{2\gamma}-\frac{\sigma}{x^{2}}\big]u_n$. Secondly, we introduce a regular weight function in the Carleman estimate for forward stochastic Grushin equations to study our inverse problem of determining two source functions simultaneously. Based on such a regular weight function, we can put the random source function $H$ on the left-hand side of this Carleman estimate, which allows us to determine $H$. However the derivatives of $H$ with respect to spatial variables still lie on the right-hand side of Carleman estimate. For this reason, the random source function $H$ to be determined could not depend on $x$ and $y$. Moreover, similar to [\ref{2012}] or [\ref{Yuan2015}], we can only determine $h$ in partial domain $I_x\times (0,T)$, since in the proof of the uniqueness result we have to differentiate the equation (\ref{1.3}) with respect to $y$, rather than $t$ as the deterministic case. This is also the result arising from the random effect of the equation.

Throughout this paper, we denote by $L^2_{\mathcal F}(0,T)$ the space of all progressively measurable stochastic process $X$ such that $\mathbb E(\int_0^T|X|^2{\rm d}t)$ $<\infty$. For a Banach space $H$, we denote by $L^2_{\mathcal F}(0,T;H)$ the Banach space consisting of all $H$-valued $\{\mathcal F_t\}_{t\geq0}$ -adapted processes $X(\cdot)$ such that $\mathbb E(| X(\cdot)|^2_{L^2(0,T;H)})$ $<\infty$, with the canonical norm; by $L_{\mathcal F}^\infty(0,T; H)$ the Banach space consisting of all $H$-valued $\{\mathcal F_t\}_{t\geq0}$-adapted bounded processes; and by $L^2_{\mathcal F}(\Omega;C([0,T]; H))$ the Banach space consisting of all $H$-valued $\{\mathcal F_t\}_{t\geq0}$-adapted continuous processes $X$ such that $\mathbb E(|X|^2_{C([0,T]; H)})<\infty$, with the canonical norm.


%

Now we state the main results in this paper. The first one is  the following null controllability for any $T$ and any degeneracy $\gamma>0$.
\begin{thm}
Let $\gamma>0$, $0\leq\sigma<\frac{1}{4}$ and $\alpha,\beta\in L_{\mathcal F}^{\infty}(0,T;L^{\infty}(I))$. Then for any $u_0\in L^2(\Omega, \mathcal F_0, \mathbb P; $ $L^2(I))$, there exists a pair $(g, G)\in L_\mathcal F^2(0, T; L^2(\omega))\times L_\mathcal F^2(0, T; L^2(I))$ such that the
corresponding solution $u$ of (\ref{1.2}) satisfies $u(T)=0$ in $I$, $\mathbb P$-a.s. for any $T>0$.
\end{thm}

{\noindent\bf Remark 1.1.}\  Condition $0\leq\sigma<\frac{1}{4}$ is used to guarantee well-posedness issues
linked to the use of the following Hardy inequality [\ref{Cannarsa2005}]
\begin{align}\label{1-1.5}
\int_0^1\frac{z^2(x)}{x^2}{\rm d}x\leq 4\int_0^1 z_x^2(x){\rm d}x,\quad \forall z\in H_0^1(0,1).
\end{align}
 Moreover, it is noted that our control domain touches the line $\{x=0\}$, which allows us to prove our controllability result for any $\gamma>0$ and any time $T>0$. However, a coming flaw with such a control domain is that the null controllability could not hold for $\sigma=\frac{1}{4}$. This is because that we need $\frac{1}{4}-\sigma>0$ to prove the Cacciopoli
inequality (\ref{1-3.45}), when our control domain $\omega$ touches the line $\{x=0\}$.

\vspace{2mm}

The other one is the following uniqueness result for our inverse source problem.
\begin{thm}
Let $\gamma>0$, $0\leq \sigma<\frac{1}{4}$,  $h\in L_\mathcal F^2(0, T; H^1(I_x))$, $H\in L^2_\mathcal F(0,T)$ and $R_1, R_2\in C^3(\overline Q_T)$ such that
\begin{align}\label{3-1.5}
&|R_i|\not=0\quad {\rm in}\ Q_T,\quad i=1,2,\\
\label{3-1.6} &\left|\nabla \left(\frac{R_2}{R_1}\right)_y\right|\leq C\left|\left(\frac{R_2}{R_1}\right)_y\right|\quad {\rm in}\ Q_T.
\end{align}
If
\begin{align}
\label{3-11.7}&u_y\big |_{\Sigma_T}=u_x\big|_{\Gamma_T}=0,\quad \mathbb P-a.s.,\\
\label{3-11.8}&u(T)=0\quad {\rm in}\  I,\quad \mathbb P-a.s.,
\end{align}
then
\begin{align}\label{3-1.7}
h(x,t)=0, \quad (x,t)\in I_x\times[0,T],\quad \mathbb P-a.s.
\end{align}
and
\begin{align}\label{3-1.8}
H(t)=0, \quad t\in [0,T],\quad \mathbb P-a.s.,
\end{align}
where $u$ is the solution of (\ref{1.3}) corresponding to $h$ and $H$.

\end{thm}

{\noindent\bf Remark 1.2.}\  Obviously, condition (\ref{3-1.6}) is correct for $\frac{R_2}{R_1}$ not depending on $y$. Or when $\Big|\nabla\ln\Big|\Big(\frac{R_2}{R_1}\Big)_y\Big|\Big|\leq C$ in $\overline Q_T$, i.e.  $\frac{R_2}{R_1}$ sufficiently smooth in $\overline Q_T$, (\ref{3-1.6}) is also correct.

\vspace{2mm}
The rest of this paper is organized as follows. In next section, we prove the well-posedness of the system \eqref{1.1}. In section 3, we establish two Carleman estimates for stochastic forward/backward Grushin equation with singular potential, respectively. In section 4, we prove the null controllability for system (\ref{1.2}), i.e. Theorem 1.1. In last section, we show the uniqueness for our inverse source problem, i.e.  Theorem 1.2.

\section{Well-posedness}
\setcounter{equation}{0}
In this section, we show the well-posedness of the following stochastic
Grushin equation with singular potential:
\bes\label{1-2.1}\left\{\begin{array}{ll}
{\rm d}u-u_{xx}{\rm d}t-x^{2\gamma}u_{yy}{\rm d}t-\frac{\sigma}{x^2}u{\rm d}t=f{\rm d}t+F{\rm d}B(t),&(x,y,t)\in Q_T,\\
u(x,y,t)=0,&(x,y,t)\in \Sigma_T,\\
u(x,y,0)=u_0(x,y),&(x,y)\in I.
\end{array}\right.
\ees In order to deal with the degeneracy and the singularity, we introduce some suitable spaces. For $\gamma>0$, we define $H_{\gamma}^1(I)$ as the completion of $C_0^{\infty}(I)$ in the norm
\begin{align*}
\|u\|_{H_{\gamma}^1(I)}=\left[\int_I\left(|u_x|^2+x^{2\gamma}|u_y|^2-\frac{\sigma}{x^2}|u|^2\right){\rm d}x{\rm d}y\right]^{\frac{1}{2}}.
\end{align*}
The Hardy inequality (\ref{1-1.5})
implies  that $H_{\gamma}^1(I)$ is a Banach
space endowed with the above norm for all $\sigma<\frac{1}{4}$.
Further we introduce
\begin{align*}
&\mathcal G_{T}=L^2_{\mathcal F}(\Omega; C([0,T]; L^2(I)))\cap L^2_{\mathcal F}(0,T;  H_{\gamma}^1(I)),\\
&\mathcal H_T=L^2_{\mathcal F}(\Omega; C([0,T]; L^2(I)))\cap L^2_{\mathcal F}(0,T; H_0^1(I)).
\end{align*}

Now, we give the definition of the weak solution of \eqref{1-2.1}.

\begin{defn}
A weak solution of \eqref{1-2.1} is a stochastic process $u\in\mathcal G_T$ such that for any $\vartheta\in C^{1}(\overline I)$, it holds that
\begin{align}
&\int_I \left[u(t)-u_0\right]\vartheta{\rm d}x{\rm d}y+\int_{Q_t} \left(u_{x}\vartheta_x+x^{2\gamma}u_{y}\vartheta_y-\frac{\sigma}{x^2}u\vartheta\right){\rm d}x{\rm d}y{\rm d}t\nonumber\\
=&\int_{Q_t} f\vartheta{\rm d}x{\rm d}y{\rm d}t+\int_{Q_t} F\vartheta{\rm d}x{\rm d}y{\rm d}B(t),\quad \mathbb P-a.s.
\end{align}
\end{defn}

\begin{thm}
Let $\gamma>0$ and $0\leq\sigma<\frac{1}{4}$. Then for any $u_0\in L^2(\Omega, \mathcal F_0, \mathbb P; L^2(I))$, system \eqref{1-2.1} admits a unique weak solution $u\in\mathcal G_T$ such that
\begin{align}\label{1-2.3}
\|u\|_{\mathcal G_T}\leq C\big(\|u_0\|_{L^2(\Omega,\mathcal F_0,\mathbb P;L^2(I))}+\|f\|_{L^2_\mathcal F(0,T;L^2(I))}+\|F\|_{L^2_\mathcal F(0,T;L^2(I))}\big),
\end{align}
where $C$ is depending on $I,T,\gamma$ and $\sigma$.
\end{thm}
{\noindent\bf Proof.}\ Letting $0<\varepsilon<1$, we consider the following approximate problem:
\bes\label{1-2.4}\left\{\begin{array}{ll}
{\rm d}u^\varepsilon-u^\varepsilon_{xx}{\rm d}t-\left(x+\varepsilon\right)^{2\gamma}u^\varepsilon_{yy}{\rm d}t-\frac{\sigma}{(x+\varepsilon)^2}u^\varepsilon{\rm d}t=f{\rm d}t+F{\rm d}B(t),&(x,y,t)\in Q_T ,\\
u^\varepsilon(x,y,t)=0,& (x,y,t)\in \Sigma_T,\\
u^\varepsilon(x,y,T)=u_0^\varepsilon(x,y),&(x,y)\in I,\end{array}\right.
\ees
where
\begin{align*}u_0^\varepsilon\rightarrow u_0\quad {\rm in}\  L^2(\Omega, \mathcal F_0, \mathbb P; L^2(I)).
\end{align*}
Then by [\ref{Peng1991SAA}], we know that  (\ref{1-2.4}) admits a unique solution $u^\varepsilon\in \mathcal H_T$ for any $0<\varepsilon<1$.

By It\^{o} formula and the equation of $u^\varepsilon$, we have
\begin{align}\label{1-2.5}
{\rm d}\big(|u^\varepsilon|^2\big)=&2u^\varepsilon{\rm d}u^\varepsilon+({\rm d}u^\varepsilon)^2\nonumber\\
=&2u^\varepsilon \left(u^\varepsilon_{xx}{\rm d}t+\left(x+\varepsilon\right)^{2\gamma}u^\varepsilon_{yy}{\rm d}t+\frac{\sigma}{(x+\varepsilon)^2}u^\varepsilon{\rm d}t+f{\rm d}t+F{\rm d}B(t)\right)+|F|^2{\rm d}t.
\end{align}
Therefore, integrating both sides of (\ref{1-2.5}) in $Q_T$ and  taking mathematical expectation in $\Omega$, we have
\begin{align}
&\mathbb E\int_I |u^\varepsilon(t)|^2{\rm d}x{\rm d}y+2\mathbb E\int_{Q_t} \left[|u^\varepsilon_{x}|^2+(x+\varepsilon)^{2\gamma}|u^\varepsilon_{y}|^2-\frac{\sigma}{(x+\varepsilon)^2}|u^\varepsilon|^2\right]{\rm d}x{\rm d}y{\rm d}t\nonumber\\
=&\mathbb E\int_I |u^\varepsilon_0|^2{\rm d}x{\rm d}y+2\mathbb E\int_{Q_t} fu^\varepsilon{\rm d}x{\rm d}y{\rm d}t+\mathbb E\int_{Q_t} |F|^2 {\rm d}x{\rm d}y{\rm d}t\nonumber\\
\leq &\mathbb E\int_I |u^\varepsilon_0|^2{\rm d}x{\rm d}y+\mathbb E\int_{Q_T} \big(|f|^2+|F|^2\big) {\rm d}x{\rm d}y{\rm d}t+\mathbb E\int_{Q_t} |u^\varepsilon|^2 {\rm d}x{\rm d}y{\rm d}t.
\end{align}
Then applying Gronwall inequality yields that
\begin{align}\label{}
&\sup_{t\in [0,T]}\mathbb E\|u^\varepsilon(t)\|^2_{L^2(I)}+\mathbb E\int_0^T\|u^\varepsilon(t)\|^2_{H^1_\gamma(I)}{\rm d}t\nonumber\\
\leq& C\mathbb E\int_{I}|u^\varepsilon_0|^2{\rm d}x+C\mathbb E\int_{Q_T}\big(|f|^2+|F|^2\big){\rm d}x{\rm d}t,
\end{align}
where $C$ is depending on $I,T$, $\gamma$ and $\sigma$, but independent of $\varepsilon$.

Similarly, we could prove for any $\varepsilon_1,\varepsilon_2\in (0,1)$ that
\begin{align}\label{}
&\sup_{t\in [0,T]}\mathbb E\|(u^{\varepsilon_1}-u^{\varepsilon_2})(t)\|^2_{L^2(I)}+\mathbb E\int_0^T\|(u^{\varepsilon_1}-u^{\varepsilon_2})(t)\|^2_{H^1_\gamma(I)}{\rm d}t\nonumber\\
\leq& C\mathbb E\int_{I}|u^{\varepsilon_1}_0-u^{\varepsilon_2}_0|^2{\rm d}x,
\end{align}
which implies that
\begin{align}\label{1-2.9}{u^\varepsilon}\rightarrow u\quad {\rm in}\ \mathcal G_T,
\end{align} due to $u_0^\varepsilon\rightarrow u_0$ in $ L^2(\Omega, \mathcal F_0, \mathbb P; L^2(I))$. Therefore by a standard limiting process we find that (\ref{1-2.1}) admits  a weak solution $u\in \mathcal G_T$ (the limit of $u^\varepsilon$ in $\mathcal G_T$) such that (\ref{1-2.3}). The uniqueness of solution could be directly deduced from (\ref{1-2.3}). \hfill$\Box$

\section{Carleman estimates for stochastic Grushin equation}
\setcounter{equation}{0}
In this section,  we will show two Carleman estimates for stochastic Grushin equation with singular potential, which will be used to study the null controllability and the inverse source problem, respectively.  One is for the backward stochastic Grushin equation with singular weight function. The other one is for the forward stochastic Grushin equation with regular weight function.

\subsection{Carleman estimate for backward stochastic Grushin equation with singular weight function}
\setcounter{equation}{0}
In this subsection, we will used a singular weight function to prove a Carlemen estimate for the backward stochastic Grushin equation with singular potential
\bes\label{3.1.1}\left\{\begin{array}{ll}
{\rm d}v+v_{xx}{\rm d}t+x^{2\gamma}v_{yy}{\rm d}t+\frac{\sigma}{x^2}v{\rm d}t=f_1{\rm d}t+F_1{\rm d}B(t),&(x,y,t)\in Q_T,\\
v(x,y,t)=0,& (x,y,t)\in \Sigma_T,\\
v(x,y,T)=v_T(x,y),&(x,y)\in I,\end{array}\right.
\ees
where $v_T\in L^2(\Omega,\mathcal F_T,\mathbb P;L^2(I))$. This Carleman estimate will be used to prove the null controllability result for (\ref{1.2}).

To formulate our Carleman estimate, we introduce some weight functions. For $\omega=(0,a)\times I_y$, we choose $\omega^{(i)}=(0,a_i)\times I_y $ for $i=1,2$ with $0<a_1<a_2<a$. Then we know that $\omega^{(1)}\subset \omega^{(2)}\subset \omega$. We define
\begin{align*}
&\phi(x,y)=e^{\lambda\psi(x,y)},\quad \varphi(x,y,t)=(e^{\lambda\psi(x,y)}-e^{2\lambda
\|\psi\|_{C(\overline I)}})\xi(t),\quad  \theta(x,y,t)=e^{s\varphi(x,y,t)},
\end{align*}
with
\begin{align}\label{3.1.4}
\psi(x,y)=x^{2+2\gamma}y(1-y)-\mu x+M, \quad \xi(t)=\frac{1}{t^4(T-t)^4}.\end{align}
Here $\mu$ is a positive constant such that
\begin{align}\label{3.1.7}
\mu >\sup_{(x,y)\in \overline I}(2+2\gamma)(x+1)^{1+2\gamma}y(1-y)+\delta_0
\end{align}
with some $\delta_0>0$, which will be specified below. $M$ is chosen sufficiently large to satisfy $\psi(x,y)>0$ for all $(x,y)\in \overline I$.
Obviously, the function $\xi$ satisfies the following essential properties
\begin{align}\label{3.1.5}\xi(t)\rightarrow +\infty\quad {\rm as}\ t\rightarrow 0^+\ {\rm or}\ T^{-}\quad {\rm and}\quad \xi>0,\quad |\xi_t|\leq C\xi^{\frac{5}{4}}.\end{align}

Our main result in the subsection is the following Carleman estimate for \eqref{3.1.1}.

\begin{thm} Let $\gamma>0$, $0\leq\sigma<\frac{1}{4}$, $v_T\in L^2(\Omega,\mathcal F_T, \mathbb P; L^2(I))$, $f_1\in L^2_{\mathcal F}(0,T; L^2(I))$, $F_1\in L^2_{\mathcal F}(0,T;L^2(I))$. Then there exist constants $\lambda_1=\lambda_1(I,T,\gamma,\sigma,\omega$, $\mu$, $M)$,
$s_1$ $=s_1(I$, $T,\gamma,\sigma, \omega,\mu,M,$ $\lambda)$ and
$C=C(I,T,\gamma,\sigma,\omega,\mu,M,\lambda)$  such that
\begin{align}\label{3.1.8}
&\mathbb{E}\int_{Q_T} s\xi\theta^2 |v_x|^2{\rm d}x{\rm
d}y{\rm d}t+\mathbb{E}\int_{Q_T}s \xi \theta^2x^{2\gamma}|v_y|^2{\rm d}x{\rm d}y{\rm d}t+\mathbb{E}\int_{Q_T} s^3\xi^3\theta^2|v|^2{\rm d}x{\rm d}y{\rm
d}t\nonumber\\
\hspace{-0.3cm}\leq & C \left[\mathbb{E}\int_{Q_T}\theta^2|f_1|^2{\rm d}x{\rm d}y{\rm d}t+\mathbb{E}\int_{Q_T}s^2\xi^2\theta^2|F_1|^2{\rm d}x{\rm d}y{\rm d}t+\mathbb{E}\int_{\omega_T}s^3\xi^3\theta^2|v|^2{\rm d}x{\rm d}y{\rm d}t\right]
\end{align}
for all $\lambda>\lambda_1$, $s>s_1$, and all $u\in \mathcal G_T$ satisfies (\ref{3.1.1}).
\end{thm}

\vspace{2mm}

Since the system (\ref{3.1.1}) is not only degenerate, but also singular on $\{x=0\}\times I_y$, we first transfer to study an approximate version of (\ref{3.1.1}). To do this, letting $0<\varepsilon<1$ and $F_1^\varepsilon\in L^2_{\mathcal F}(0,T;H_0^1(I))$, $v^\varepsilon_T\in L^2(\Omega,\mathcal F_T,\mathbb P;H_0^1(I))$ such that
\begin{align*}
&F_1^\varepsilon\rightarrow F_1\quad {\rm in }\  L^2_{\mathcal F}(0,T;L^2(I)),\\
&v^\varepsilon_T\rightarrow v_T\quad {\rm in}\ L^2(\Omega,\mathcal F_T,\mathbb P;L^2(I)),
\end{align*}
we then consider
\bes\label{3.1.2}\left\{\begin{array}{ll}
{\rm d}v^\varepsilon+v^\varepsilon_{xx}{\rm d}t+\left(x+\varepsilon\right)^{2\gamma}v^\varepsilon_{yy}{\rm d}t+\frac{\sigma}{(x+\varepsilon)^2}v^\varepsilon{\rm d}t=f_1{\rm d}t+F_1^\varepsilon{\rm d}B(t),&(x,y,t)\in Q_T ,\\
v^\varepsilon(x,y,t)=0,& (x,y,t)\in \Sigma_T,\\
v^\varepsilon(x,y,T)=v_T^\varepsilon(x,y),&(x,y)\in I.\end{array}\right.
\ees
 According to the standard theory for stochastic parabolic equation, e.g. [\ref{Zhang2009},\ref{Zhou1992}], we know that the system (\ref{3.1.2}) admits a unique solution $v^{\varepsilon}\in\mathcal H_T$.
Set
\begin{align*}
\widehat\varphi(x,y,t)=\varphi(x+\varepsilon,y,t),\quad \widehat\theta(x,y,t)=\theta(x+\varepsilon,y,t).
\end{align*}
In the sequel, $\widehat \phi$ and $\widehat \psi$ are defined analogously. Then we have the following weighted identity for (\ref{3.1.2}).

\begin{lem}Let $\tau$ be a constant such that $2<\tau<3$. Assume that $v^\varepsilon$ is an $H^2(\mathbb R^2)$-valued continuous semimartingale. Set $l=s\widehat\varphi$, $z=\widehat \theta v^\varepsilon$ and
\begin{align*}
{{P}}_1=&{\rm d}z-2 l_xz_x{\rm d}t-2\left(x+\varepsilon\right)^{2\gamma} l_yz_y{\rm d}t-\tau l_{xx}z{\rm d}t,\\
{{P}}_2=&z_{xx}+\left(x+\varepsilon\right)^{2\gamma}z_{yy}+ l_x^2z+\left(x+\varepsilon\right)^{2\gamma}l_y^2z+\frac{\sigma}{(x+\varepsilon)^2}z,\nonumber\\
{{P}}=&(\tau-1) l_{xx}z- l_tz-\left(x+\varepsilon\right)^{2\gamma} l_{yy}z.
\end{align*}
Then for a.e. $(x,y)\in \mathbb R^2$, it holds that
\begin{align}\label{3.1.17}
&{ P}_2 \widehat\theta\left[{\rm d}v^\varepsilon+v^\varepsilon_{xx}{\rm d}t+\left(x+\varepsilon\right)^{2\gamma}v^\varepsilon_{yy}{\rm d}t+\frac{\sigma}{(x+\varepsilon)^2}v^\varepsilon{\rm d}t\right]\nonumber\\
=& |{P}_2| ^2{\rm d}t+{P}_2{P}{\rm d}t+\sum_{i=1}^{5}X_i{\rm d}t+{\rm d}Y+\{\cdot\}_x+\{\cdot\cdot\}_y+J,\quad \mathbb P-a.s.,
\end{align}
where
\begin{align*}
X_1=&\left[(\tau+1)l_{xx}-(x+\varepsilon)^{2\gamma}l_{yy}\right]z_x^2, \nonumber\\
X_2=&\left[-2\gamma(x+\varepsilon)^{2\gamma-1}l_{x}+(\tau-1)(x+\varepsilon)^{2\gamma}l_{xx}+(x+\varepsilon)^{4\gamma}l_{yy}\right]z_y^2,\nonumber\\
X_3=&4\left[\gamma(x+\varepsilon)^{2\gamma-1}l_y+(x+\varepsilon)^{2\gamma}l_{xy}\right]z_x z_y,\nonumber\\
X_4=&\left[(3-\tau)l^2_{x}l_{xx}+2\gamma(x+\varepsilon)^{2\gamma-1}l_x l_y^2+3(x+\varepsilon)^{4\gamma} l_y^2 l_{yy}\right]z^2\nonumber\\
&+(x+\varepsilon)^{2\gamma}\left[4 l_x l_y l_{xy}+ l_{x}^2 l_{yy}+
(1-\tau)l_{xx} l_y^2\right]z^2\nonumber\\
&+\left[(1-\tau)\frac{\sigma}{(x+\varepsilon)^2} l_{xx}-\frac{2\sigma}{(x+\varepsilon)^3}l_x
+\frac{\sigma}{(x+\varepsilon)^{2-2\gamma}}l_{yy}\right]z^2,\nonumber\\
X_5=&\left[-l_x l_{xt}-(x+\varepsilon)^{2\gamma} l_y l_{yt}
-\frac{\tau}{2} l_{xxxx}-\frac{\tau}{2}(x+\varepsilon)^{2\gamma} l_{xxyy}\right]z^2,\nonumber\\
Y=&-\frac{1}{2}z_x^2-\frac{1}{2}(x+\varepsilon)^{2\gamma}z_y^2+\frac{1}{2}\left[ l_x^2+(x+\varepsilon)^{2\gamma} l_y^2+\frac{\sigma}{(x+\varepsilon)^2}\right]z^2,\nonumber\\
\{\cdot\}=&z_x {\rm d}z+\bigg[-l_x z_x^2+(x+\varepsilon)^{2\gamma} l_x z^2_y- l_x^3z^2-(x+\varepsilon)^{2\gamma} l_x l_y^2z^2\nonumber\\
&-\frac{\sigma}{(x+\varepsilon)^2} l_x z^2-2(x+\varepsilon)^{2\gamma} l_y z_x z_y-\tau  l_{xx}z z_x+\frac{\tau}{2} l_{xxx}z^2\bigg]{\rm d}t,\nonumber\\
\{\cdot\cdot\}=&(x+\varepsilon)^{2\gamma}z_y{\rm d}z+\Big[-2(x+\varepsilon)^{2\gamma} l_x z_x z_y+(x+\varepsilon)^{2\gamma} l_{y}z_x^2\nonumber\\
&-(x+\varepsilon)^{4\gamma} l_{y}z_y^2-(x+\varepsilon)^{2\gamma} l_{x}^2 l_yz^2-(x+\varepsilon)^{4\gamma}l_y^3z^2\nonumber\\
&-\frac{\sigma}{(x+\varepsilon)^{2-2\gamma}} l_yz^2-\tau(x+\varepsilon)^{2\gamma}l_{xx}zz_{y}+\frac{\tau}{2}(x+\varepsilon)^{2\gamma}l_{xxy}z^2\Big]{\rm d}t,\\
J=&\frac{1}{2}({\rm d}z_x)^2+\frac{1}{2}(x+\varepsilon)^{2\gamma}({\rm d}z_y)^2-\frac{1}{2}\left[l_x^2+(x+\varepsilon)^{2\gamma}l_y^2+\frac{\sigma}{(x+\varepsilon)^2}\right]({\rm d}z)^2.
\end{align*}
\end{lem}

{\noindent\bf Proof.}\  Notice that $\widehat\theta=e^{l}$, $l=s\widehat\varphi$ and $z=\widehat\theta v^\varepsilon$. Then we have
\begin{align*}
\widehat\theta\left[{\rm d}v^\varepsilon+v^\varepsilon_{xx}{\rm d}t+\left(x+\varepsilon\right)^{2\gamma}v^\varepsilon_{yy}{\rm d}t+\frac{\sigma}{(x+\varepsilon)^2}v^\varepsilon{\rm d}t\right]= P_1+(P_2+P){\rm d}t.
\end{align*}
Hence
\begin{align}\label{3.1.11}
P_2\widehat\theta \left[{\rm d}v^\varepsilon+v^\varepsilon_{xx}{\rm d}t+\left(x+\varepsilon\right)^{2\gamma}v^\varepsilon_{yy}{\rm d}t+\frac{\sigma}{(x+\varepsilon)^2}v^\varepsilon{\rm d}t\right]= P_1  P_2+| P_2|^2{\rm d}t+ P_2 P{\rm d}t.
\end{align}

We easily see that
\begin{align}\label{3.1.12}
P_1P_2=P_2 {\rm d}z-2l_x z_x P_2{\rm d}t-2(x+\varepsilon)^{2\gamma} l_y z_y P_2{\rm d}t-\tau  l_{xx} z P_2{\rm d}t.
\end{align}
Now we calculate the terms on the right-hand side of \eqref{3.1.12} one by one. For the first one, by It\^{o}'s formula, we have
\begin{align}\label{3.1.13}
P_2{\rm d}z=&\left[z_{xx}+\left(x+\varepsilon\right)^{2\gamma}z_{yy}+ l_x^2z+\left(x+\varepsilon\right)^{2\gamma}l_y^2z+\frac{\sigma}{(x+\varepsilon)^2}z\right]{\rm d}z\nonumber\\
=&(z_{x}{\rm d}z)_x-\frac{1}{2}{\rm d}(z_x^2)+\frac{1}{2}({\rm d}z_x)^2+\left[(x+\varepsilon)^{2\gamma}z_y{\rm d}z\right]_y-\frac{1}{2}{\rm d}[(x+\varepsilon)^{2\gamma}z_y^2]\nonumber\\
&+\frac{1}{2}(x+\varepsilon)^{2\gamma}({\rm d}z_y)^2+\frac{1}{2}{\rm d}\left(l_x^2z^2\right)-l_xl_{xt}z^2{\rm d}t -\frac{1}{2}l_x^2({\rm d}z)^2\nonumber\\
&+\frac{1}{2}{\rm d}\left[(x+\varepsilon)^{2\gamma}\ l_y^2z^2\right]-(x+\varepsilon)^{2\gamma}l_y l_{yt}z^2{\rm d}t-\frac{1}{2}(x+\varepsilon)^{2\gamma}l_y^2({\rm d}z)^2\nonumber\\
&+\frac{1}{2}{\rm d}\left[\frac{\sigma}{(x+\varepsilon)^2}z^2\right]-\frac{1}{2}\frac{\sigma}{(x+\varepsilon)^2}({\rm d}z)^2.
\end{align}
By a direct calculation, we have
\begin{align}\label{3.1.14}
&-2l_xz_x P_2{\rm d}t\nonumber\\
=&-2l_xz_x\left[z_{xx}+\left(x+\varepsilon\right)^{2\gamma}z_{yy}+ l_x^2z+\left(x+\varepsilon\right)^{2\gamma}l_y^2z+\frac{\sigma}{(x+\varepsilon)^2}z\right]{\rm d}t\nonumber\\
=&-(l_xz_x^2)_x{\rm d}t+l_{xx}z_x^2{\rm d}t
-2\left[(x+\varepsilon)^{2\gamma}\ l_xz_xz_y\right]_y{\rm d}t+\left[(x+\varepsilon)^{2\gamma}\ l_xz_y^2\right]_x{\rm d}t\nonumber\\&
+2(x+\varepsilon)^{2\gamma}l_{xy}z_xz_y{\rm d}t-\left[2\gamma (x+\varepsilon)^{2\gamma-1}l_x+(x+\varepsilon)^{2\gamma}l_{xx}\right]z_y^2{\rm d}t-\left(l^3_xz^2\right)_x{\rm d}t\nonumber\\
&+3l_x^2l_{xx}z^2{\rm d}t-\left[(x+\varepsilon)^{2\gamma}l_x l_y^2z^2\right]_x{\rm d}t+2(x+\varepsilon)^{2\gamma}l_xl_yl_{xy}z^2{\rm d}t\nonumber\\
&+\left[2\gamma(x+\varepsilon)^{2\gamma-1}l_xl_y^2
+(x+\varepsilon)^{2\gamma}l_{xx}l_y^2\right]z^2{\rm d}t-\left[\frac{\sigma}{(x+\varepsilon)^2}l_xz^2\right]_x{\rm d}t\nonumber\\
&+\left[
\frac{\sigma}{(x+\varepsilon)^2}l_{xx}-\frac{2\sigma}{(x+\varepsilon)^3}l_x\right]z^2{\rm d}t
\end{align}
and
\begin{align}\label{3.1.15}
&-2(x+\varepsilon)^{2\gamma}l_yz_y P_2{\rm d}t\nonumber\\
=&-2(x+\varepsilon)^{2\gamma}l_yz_y\left[z_{xx}+\left(x+\varepsilon\right)^{2\gamma}z_{yy}+ l_x^2z+\left(x+\varepsilon\right)^{2\gamma} l_y^2z+\frac{\sigma}{(x+\varepsilon)^2}z\right]{\rm d}t\nonumber\\
=&-2\left[\left(x+\varepsilon\right)^{2\gamma}l_yz_xz_y\right]_x{\rm d}t+\left[\left(x+\varepsilon\right)^{2\gamma}l_yz_x^2\right]_y{\rm d}t\nonumber\\
&+\left[4\gamma (x+\varepsilon)^{2\gamma-1}l_y+2(x+\varepsilon)^{2\gamma}l_{xy}\right]z_xz_y{\rm d}t-(x+\varepsilon)^{2\gamma}l_{yy}z_x^2{\rm d}t\nonumber\\
&-\left[(x+\varepsilon)^{4\gamma}l_yz_y^2\right]_y{\rm d}t+(x+\varepsilon)^{4\gamma}l_{yy}z_y^2{\rm d}t-\left[(x+\varepsilon)^{2\gamma}l_x^2l_yz^2\right]_y{\rm d}t\nonumber\\
&+(x+\varepsilon)^{2\gamma}\left(l_x^2l_y\right)_yz^2{\rm d}t-\left[(x+\varepsilon)^{4\gamma}l_y^3z^2\right]_y{\rm d}t+3(x+\varepsilon)^{4\gamma}l_y^2l_{yy}z^2{\rm d}t\nonumber\\
&-\left[\frac{\sigma}{(x+\varepsilon)^{2-2\gamma}}l_y z^2\right]_y{\rm d}t+\frac{\sigma}{(x+\varepsilon)^{2-2\gamma}}l_{yy}z^2{\rm d}t.
\end{align}
The last term can be rewritten as
\begin{align}\label{3.1.16}
&-\tau l_{xx}z P_2{\rm d}t\nonumber\\
=&-\tau l_{xx}z\left[z_{xx}+\left(x+\varepsilon\right)^{2\gamma}z_{yy}+ l_x^2z+\left(x+\varepsilon\right)^{2\gamma}l_y^2z+\frac{\sigma}{(x+\varepsilon)^2}z\right]{\rm d}t\nonumber\\
=&-\tau\left(l_{xx}zz_x\right)_x{\rm d}t+\frac{\tau}{2}\left(l_{xxx}z^2\right)_x{\rm d}t-\frac{\tau}{2}l_{xxxx}z^2{\rm d}t+\tau l_{xx}z_x^2{\rm d}t-\tau\left[(x+\varepsilon)^{2\gamma}l_{xx}zz_y\right]_y{\rm d}t\nonumber\\
&+\frac{\tau}{2}\left[(x+\varepsilon)^{2\gamma}l_{xxy}z^2\right]_y{\rm d}t-\frac{\tau}{2}(x+\varepsilon)^{2\gamma}l_{xxyy}z^2{\rm d}t+\tau(x+\varepsilon)^{2\gamma}l_{xx}z_y^2{\rm d}t\nonumber\\
&-\tau l^2_xl_{xx}z^2-\tau(x+\varepsilon)^{2\gamma}l_{xx} l_y^2z^2-\tau\frac{\sigma}{(x+\varepsilon)^2} l_{xx} z^2{\rm d}t.
\end{align}
Combining (\ref{3.1.11})-(\ref{3.1.16}), we can obtain (\ref{3.1.17}) and then complete the proof of Lemma 3.2. \hfill $\Box$

\vspace{2mm}

Now, integrating both sides of (\ref{3.1.17}) in $Q_T$, taking mathematical expectation in $\Omega$ and using $\widehat \theta(x,y,0)=\widehat \theta(x,y,T)=0$ in $I$, we obtain that
\begin{align}\label{3.1.18}
&\mathbb{E}\int_{Q_T} P_2 \widehat\theta\left[{\rm d}v^\varepsilon+v^\varepsilon_{xx}{\rm d}t+\left(x+\varepsilon\right)^{2\gamma}v^\varepsilon_{yy}{\rm d}t+\frac{\sigma}{(x+\varepsilon)^2}v^\varepsilon{\rm d}t\right]{\rm d}x{\rm d}y\nonumber\\
\geq &\frac{1}{2}\mathbb{E}\int_{Q_T}|P_2|^2{\rm d}x{\rm d}y{\rm d}t-\frac{1}{2}\mathbb{E}\int_{Q_T}|P|^2{\rm d}x{\rm d}y{\rm d}t+\sum_{i=1}^{5}\mathbb{E}\int_{Q_T}X_i{\rm d}x{\rm d}y{\rm d}t\nonumber\\
&+\mathbb{E}\int_{Q_T}\left(\{\cdot\}_x+\{\cdot\cdot\}_y\right){\rm d}x{\rm d}y+\mathbb{E}\int_{Q_T}J{\rm d}x{\rm d}y.
\end{align}
In the following we estimate the last three terms in (\ref{3.1.18}) one by one.

\begin{lem} There exists constant $C=C(I,T,\gamma,\sigma,\mu,M)$  such that
\begin{align}\label{3.15}
\sum_{i=1}^{5}\mathbb{E}\int_{Q_T}X_i{\rm d}x{\rm d}y{\rm d}t\geq &C\mathbb{E}\int_{Q_T}s^3\lambda^4\widehat\phi^3 \xi^3|z|^2{\rm d}x{\rm d}y{\rm d}t+C\mathbb{E}\int_{Q_T} s\lambda^2\widehat\phi\xi |z_x|^2{\rm d}x{\rm
d}y{\rm d}t\nonumber\\
&+C\mathbb{E}\int_{Q_T}s\lambda^2(x+\varepsilon)^{2\gamma}\widehat\phi\xi |z_y|^2{\rm d}x{\rm d}y{\rm d}t
\end{align}
for all large $\lambda$ and $s$.
\end{lem}

\noindent{\bf Proof.} Notice that $\widehat\psi(x,y)=(x+\varepsilon)^{2+2\gamma}y(1-y)-\mu(x+\varepsilon)+M$. Together with (\ref{3.1.7}), we obtain the
following properties of $\widehat\psi$:
\begin{align}\label{3.1.19}\left\{\begin{array}{l}
\widehat\psi_x <-\delta_0,\quad \widehat\psi_x\widehat\psi_{xxx}\leq 0,\\
|\widehat\psi_{xxxx}|\leq C(x+\varepsilon)^{-2},\quad |\widehat\psi_{y}|+|\widehat\psi_{yy}|\leq C(x+\varepsilon)^{2+2\gamma},\\
|\widehat\psi_{xy}|+|\widehat\psi_{xx}|+|\widehat\psi_{xxy}|+|\widehat\psi_{xyy}|+|\widehat\psi_{xxyy}|\leq C(x+\varepsilon)^{2\gamma}.
\end{array}\right.\end{align}
Recalling $l=s\widehat \varphi$, we have
\begin{align}\label{3.17}\left\{\begin{array}{l}
X_1=(\tau+1)s\lambda^2\widehat\phi\widehat\psi_x^2 \xi|z_x|^2+K_1|z_x|^2,\\
X_2=\left[s\lambda^2(x+\varepsilon)^{4\gamma}\widehat\psi_y^2+(\tau-1)s\lambda^2(x+\varepsilon)^{2\gamma}\widehat\psi_x^2-2\gamma s\lambda(x+\varepsilon)^{2\gamma-1}\widehat \psi_x\right]\widehat\phi \xi |z_y|^2\\
\quad\quad\ \ +K_2|z_y|^2,\\
X_3=4s\lambda^2(x+\varepsilon)^{2\gamma}\widehat\phi\widehat\psi_x\widehat\psi_y \xi z_xz_y+K_3 z_xz_y,
\end{array}
\right.
\end{align}
where
\begin{align*}
K_1=&\left[-s\lambda^2(x+\varepsilon)^{2\gamma}\widehat\psi_y^2+s\lambda\left((\tau+1)\widehat\psi_{xx}-(x+\varepsilon)^{2\gamma}\widehat\psi_{yy}\right)\right]\widehat\phi \xi,\nonumber\\
K_2=&s\lambda \left[(\tau-1)(x+\varepsilon)^{2\gamma}\widehat\psi_{xx}+(x+\varepsilon)^{4\gamma}\widehat\psi_{yy}\right]\widehat \phi \xi ,\nonumber\\
K_3=&4s\lambda\left[(x+\varepsilon)^{2\gamma}\widehat\psi_{xy}+\gamma(x+\varepsilon)^{2\gamma-1}
\widehat\psi_y\right]\widehat\phi\xi,
\end{align*}
satisfy
\begin{align}\left\{\begin{array}{l}
|K_1|\leq Cs\lambda^2\widehat\phi\xi,\\
|K_2|\leq Cs\lambda(x+\varepsilon)^{2\gamma}\widehat\phi\xi,\\
|K_3|\leq Cs\lambda(x+\varepsilon)^{2\gamma}\widehat\phi\xi.\end{array}\right.
\end{align}
 due to (\ref{3.1.19}). By Young's inequality, we obtain for all $\epsilon>0$ that
\begin{align}\label{1-3.19}
\left|4s\lambda^2(x+\varepsilon)^{2\gamma}\widehat\phi\widehat\psi_x\widehat\psi_y \xi z_xz_y\right|
\leq\epsilon s\lambda^2\widehat\psi_x^2 \widehat\phi\xi|z_x|^2+C(\epsilon)s\lambda^2(x+\varepsilon)^{4\gamma}\widehat\psi_y^2 \widehat\phi\xi |z_y|^2.
\end{align}
Therefore, by (\ref{3.17})-(\ref{1-3.19}) we have the following estimate
\begin{align}\label{3.1.22}
&\sum_{i=1}^3\mathbb{E}\int_{Q_T} X_i{\rm d}x{\rm d}y{\rm d}t\nonumber\\
\geq&\mathbb{E}\int_{Q_T}\left[(\tau+1-\epsilon) s\lambda^2\widehat\psi_x^2-Cs\lambda^2-Cs\lambda\right]\widehat\phi \xi|z_x|^2{\rm d}x{\rm d}y{\rm d}t+\mathbb{E}\int_{Q_T}\left[
(\tau-1)s\lambda^2(x+\varepsilon)^{2\gamma}\widehat\psi_x^2\right.\nonumber\\
&\left.-2\gamma s\lambda(x+\varepsilon)^{2\gamma-1}\widehat \psi_x-C(\epsilon)s\lambda^2(x+\varepsilon)^{4\gamma}-Cs\lambda(x+\varepsilon)^{2\gamma}\right]
\widehat\phi\xi |z_y|^2{\rm d}x{\rm d}y{\rm d}t.
\end{align}Fixing $0<\epsilon<\frac{1}{2}$ and choosing $\delta_0$ sufficiently large to
satisfy
\begin{align}\label{3.1.23}\left\{\begin{array}{l}
\big(\frac{1}{2}-\epsilon\big)\delta_0^2s\lambda^2-Cs\lambda^2-Cs\lambda>0,\\
(\tau-2)\delta_0^2s\lambda^2-C(\epsilon)2^{2\gamma}s\lambda^2-Cs\lambda>0,
\end{array}\right.
\end{align}
and noticing that $\widehat\psi_x<0$, we further find that
\begin{align}\label{3.1.24}
&\sum_{i=1}^3\mathbb{E}\int_{Q_T} X_i{\rm d}x{\rm d}y{\rm d}t\nonumber\\
\geq &\left(\tau+\frac{1}{2}\right)\mathbb{E}\int_{Q_T} s\lambda^2\widehat\psi_x^2 \widehat\phi\xi|z_x|^2{\rm d}x{\rm
d}y{\rm d}t+\mathbb{E}\int_{Q_T}s\lambda^2(x+\varepsilon)^{2\gamma}\widehat\psi_x^2 \widehat\phi\xi |z_y|^2{\rm d}x{\rm d}y{\rm d}t.
\end{align}
 By definitions of $l$, $\widehat \varphi$, we have the following estimate for $X_4$:
\begin{align}
X_4=&s^3\lambda^4\left[(3-\tau)\widehat\psi_x^4+3(x+\varepsilon)^{4\gamma}\widehat\psi_y^4+(6-\tau)(x+\varepsilon)^{2\gamma}\widehat\psi_x^2\widehat\psi_y^2\right]\widehat\phi^3\xi^3|z|^2
\nonumber\\
&+\frac{s\sigma}{(x+\varepsilon)^2}\left[(1-\tau)\left(\lambda^2\widehat\psi_x^2+\lambda\widehat\psi_{xx}\right)-\frac{2}{(x+\varepsilon)}
\lambda\widehat\psi_x\right]\widehat\phi \xi |z|^2+K_4|z|^2,
\end{align}
where
\begin{align*}
K_4=&s^3\lambda^3\Big[(3-\tau)\widehat\psi_x^2\widehat\psi_{xx}+2\gamma(x+\varepsilon)^{2\gamma-1}\widehat\psi_x\widehat\psi_y^2+3(x+\varepsilon)^{4\gamma}\widehat\psi_y^2\widehat\psi_{yy}\\
&+(x+\varepsilon)^{2\gamma}\Big(4\widehat\psi_x\widehat\psi_y\widehat\psi_{xy}+\widehat\psi_x^2\widehat\psi_{yy}+(1-\tau)\widehat\psi_y^2\widehat\psi_{xx}\Big)\Big]\widehat\phi^3\xi^3\\
&+s\sigma(x+\varepsilon)^{-2+2\gamma}
\left(\lambda^2\widehat\psi_y^2+\lambda\widehat\psi_{yy}\right)
\widehat\phi \xi
\end{align*}
satisfies
\begin{align}
|K_4|\leq Cs^3\lambda^3\widehat\phi^3\xi^3.
\end{align}
Then we obtain that
\begin{align}\label{3.22}
&\mathbb{E}\int_{Q_T}X_4{\rm d}x{\rm d}y{\rm d}t\nonumber\\
\geq &\mathbb E\int_{Q_T}\left[(3-\tau)s^3\lambda^4\widehat\psi_x^4\widehat\phi^3\xi^3-\tau
\frac{\sigma}{(x+\varepsilon)^2}s\lambda^2 \widehat\psi_x^2\widehat\phi\xi -Cs^3\lambda^3\widehat\phi^3\xi^3\right]|z|^2{\rm d}x{\rm d}y{\rm d}t.
\end{align}
Moreover, by (\ref{3.1.5}) and (\ref{3.1.19}) we have
\begin{align}\label{3.23}
&\mathbb{E}\int_{Q_T}X_5{\rm d}x{\rm d}y{\rm d}t\nonumber\\
=&-\mathbb{E}\int_{Q_T}s^2\lambda^2\left[\widehat\psi_x^2+(x+\varepsilon)^{2\gamma}\widehat\psi_y^2\right]\widehat\phi^2\xi\xi_t|z|^2{\rm d}x{\rm d}y{\rm d}t\nonumber\\
&-\frac{\tau}{2}\mathbb{E}\int_{Q_T} s\left[\lambda^4\widehat\psi_x^4+6\lambda^3\widehat\psi_x^2\widehat\psi_{xx}+\lambda^2(3\widehat\psi_{xx}^2+4\widehat\psi_x\widehat\psi_{xxx})+\lambda\widehat\psi_{xxxx}\right]
\widehat\phi\xi|z|^2{\rm d}x{\rm dy}{\rm d}t\nonumber\\
&-\frac{\tau}{2}\mathbb{E}\int_{Q_T}(x+\varepsilon)^{2\gamma}s\left[\lambda^4\widehat\psi_x^2\widehat\psi_y^2+\lambda^3(4\widehat\psi_x\widehat\psi_{y}\widehat\psi_{xy}+\widehat\psi_{xx}\widehat\psi_y^2
+\widehat\psi_x^2\widehat\psi_{yy})\right.\nonumber\\
&\left.\quad\quad\ \quad\quad+\lambda^2(2\widehat\psi_{xy}^2+2\widehat\psi_y\widehat\psi_{xxy}+2\widehat\psi_x\widehat\psi_{xyy}+\widehat\psi_{xx}\widehat\psi_{yy})
+\lambda\widehat\psi_{xxyy}\right]
\widehat\phi\xi|z|^2{\rm d}x{\rm dy}{\rm d}t\nonumber\\
\geq &-C(\delta_0)\mathbb{E}\int_{Q_T}\left(s^2\lambda^2\widehat\phi^2\xi^3+s\lambda^4\widehat\phi\xi\right)|z|^2{\rm d}x{\rm d}y{\rm d}t\nonumber\\
&-C\mathbb E\int_{Q_T}s\lambda\frac{1}{(x+\varepsilon)^{2}}\widehat\phi\xi|z|^2{\rm d}x{\rm d}y{\rm d}t.
\end{align}
By Hardy inequality (\ref{1-1.5}), we have
\begin{align}\label{1-3.25}
&-\mathbb E\int_{Q_T}s\lambda\frac{1}{(x+\varepsilon)^{2}}\widehat\phi\xi|z|^2{\rm d}x{\rm d}y{\rm d}t\geq -4\mathbb E\int_{Q_T}s\lambda \xi\left|\left(\widehat\phi^{\frac{1}{2}} z\right)_x\right|^2{\rm d}x{\rm d}y{\rm d}t\nonumber\\
\geq&-C\mathbb E\int_{Q_T}s\lambda\widehat\phi\xi|z_x|^2{\rm d}x{\rm d}y{\rm d}t-C\mathbb E\int_{Q_T}s\lambda^3\widehat\psi_x^2\widehat\phi\xi|z|^2{\rm d}x{\rm d}y{\rm d}t.
\end{align}
Then, it follows from (\ref{3.1.24}), (\ref{3.22})-(\ref{1-3.25})  that
\begin{align}\label{3.1.26}
&\sum_{i=1}^{5}\mathbb{E}\int_{Q_T}X_i{\rm d}x{\rm d}y{\rm d}t\nonumber\\
\geq &\mathbb{E}\int_{Q_T}\left[(3-\tau)\delta_0^4s^3\lambda^4-C(\delta_0)s^3\lambda^3-C(\delta_0)s\lambda^4\right]\widehat\phi^3 \xi^3|z|^2{\rm d}x{\rm d}y{\rm d}t\nonumber\\
&+\tau\mathbb{E}\int_{Q_T} s\lambda^2\left[|z_x|^2-\frac{\sigma}{(x+\varepsilon)^2}|z|^2\right]\widehat\psi_x^2 \widehat\phi\xi{\rm d}x{\rm
d}y{\rm d}t+\frac{1}{2}\mathbb{E}\int_{Q_T} \left(s\lambda^2\widehat\psi_x^2-Cs\lambda\right)\widehat\phi\xi|z_x|^2{\rm d}x{\rm d}y{\rm d}t\nonumber\\
&+\mathbb{E}\int_{Q_T}s\lambda^2(x+\varepsilon)^{2\gamma}\widehat\psi_x^2 \widehat\phi\xi |z_y|^2{\rm d}x{\rm d}y{\rm d}t.
\end{align}
Then noticing that $\tau<3$ and choosing $\lambda$ and $s$ sufficiently large, we could obtain the desired estimate (\ref{3.15}). \hfill$\Box$

\begin{lem}  There exists constant $C=C(I,T,\gamma,\sigma,\mu,M)$  such that
\begin{align}\label{1-3.27}
&\mathbb{E}\int_{Q_T}\left(\{\cdot\}_x+\{\cdot\cdot\}_y\right){\rm d}x{\rm d}y\geq -C \mathbb E\int_{\Gamma_T} s\lambda \widehat \phi\xi |z_x|^2{\rm d}y{\rm d}t.
\end{align}
\end{lem}

\noindent{\bf Proof.} From the homogeneous Dirichlet boundary condition in (\ref{3.1.2}), it follows that
\begin{align}\label{3.1.38}\left\{\begin{array}{ll}
z_x(x,0,t)=z_x(x,1,t)=0,&(x,t)\in I_x\times (0,T),\\
z_y(0,y,t)=z_y(1,y,t)=0,&(y,t)\in I_y\times (0,T),\\
z_t(x,y,t)=0,& (x,y,t)\in \Sigma_T.\end{array}\right.
\end{align}
Moreover, we easily see  that
\begin{align}\label{3.1.39}
\left\{\begin{array}{ll}\widehat\varphi_x(0,y,t)\leq 0,\quad \widehat\varphi_x(1,y,t)\leq 0,&(y,t)\in I_y\times (0,T),\\
\widehat\varphi_y(x,0,t)\geq 0,\quad \widehat\varphi_y(x,1,t)\leq  0,&(x,t)\in I_x\times (0,T).\end{array}\right.
\end{align}
Therefore, by (\ref{3.1.38}) we have
\begin{align}\label{3.1.41}
&\mathbb{E}\int_{Q_T}\left(\{\cdot\}_x+\{\cdot\cdot\}_y\right){\rm d}x{\rm d}y\nonumber\\
=&-\mathbb{E}\int_0^T\int_{I_y}\left[s\widehat\varphi_xz_x^2\right]_{x=0}^{x=1}{\rm
d}y{\rm d}t
-\mathbb{E}\int_0^T\int_{I_x}\left[s(x+\varepsilon)^{4\gamma}\widehat\varphi_y z_y^2\right]_{y=0}^{y=1}{\rm
d}x{\rm d}t.
\end{align}
Finally, by (\ref{3.1.39}) and (\ref{3.1.41}) we obtain (\ref{1-3.27}). \hfill$\Box$

\begin{lem} There exists constant $C=C(I,T,\gamma,\sigma,\mu,M)$  such that
\begin{align}\label{3.1.30}
\mathbb{E}\int_{Q_T}J{\rm d}x{\rm d}y\geq-C\mathbb{E}\int_{Q_T}s^2\lambda^2\widehat\phi^2\xi^2\widehat\theta^2|F^\varepsilon_1|^2{\rm d}x{\rm d}y{\rm d}t.
\end{align}
\end{lem}

{\noindent\bf Proof.}\ By using $({\rm d}z)^2=\widehat\theta^2|F^\varepsilon_1|^2{\rm d}t$ and Hardy inequality, we find that
\begin{align}\label{3.1.28}
&\mathbb{E}\int_{Q_T}\frac{\sigma}{(x+\varepsilon)^2}({\rm d}z)^2{\rm d}x{\rm d}y\leq 4\sigma \mathbb{E}\int_{Q_T}\left|\left(\widehat\theta F_1^\varepsilon\right)_x\right|^2{\rm d}x{\rm d}y{\rm d}t\nonumber\\
\leq&4\sigma \mathbb{E}\int_{Q_T}\left( l_x^2\widehat\theta^2|F^\varepsilon_1|^2+2 l_x \widehat\theta^2F^\varepsilon_1 F^\varepsilon_{1,x}+\widehat\theta^2|F^\varepsilon_{1,x}|^2\right){\rm d}x{\rm d}y{\rm d}t.
\end{align}
By (\ref{3.1.28}) and \begin{align}({\rm d}z_x)^2=l_x^2\widehat\theta^2|F^\varepsilon_1|^2{\rm d}t+2l_x\widehat\theta^2F^\varepsilon_1F^\varepsilon_{1,x}{\rm d}t
+\widehat\theta^2|F^\varepsilon_{1,x}|^2{\rm d}t,\end{align} we further obtain
\begin{align*}
&\mathbb{E}\int_{Q_T}J{\rm d}x{\rm d}y\nonumber\\
\geq& \mathbb{E}\int_{Q_T} l_x\widehat\theta^2F^\varepsilon_1 F^\varepsilon_{1,x}{\rm d}x{\rm d}y{\rm d}t
+\frac{1}{2}\mathbb{E}\int_{Q_T}\widehat\theta^2|F^\varepsilon_{1,x}|^2{\rm d}x{\rm d}y{\rm d}t
-\frac{1}{2}\mathbb{E}\int_{Q_T}(x+\varepsilon)^{2\gamma} l_y^2\widehat\theta^2 |F^\varepsilon_1|^2{\rm d}x{\rm d}y{\rm d}t\nonumber\\
&-2\sigma \mathbb{E}\int_{Q_T}\left(l_x^2\widehat\theta^2 |F^\varepsilon_1|^2+2 l_x\widehat\theta^2 F^\varepsilon_1 F^\varepsilon_{1,x}+\widehat\theta^2|F^\varepsilon_{1,x}|^2\right){\rm d}x{\rm d}y{\rm d}t
\end{align*}
\begin{align}\label{3.1.29}
\geq&-\epsilon \mathbb{E}\int_{Q_T}\widehat\theta^2|F^\varepsilon_{1,x}|^2{\rm d}x{\rm d}y{\rm d}t-C(\epsilon) \mathbb{E}\int_{Q_T}l_x^2\widehat\theta^2|F^\varepsilon_1|^2{\rm d}x{\rm d}y{\rm d}t
+\frac{1}{2}\mathbb{E}\int_{Q_T}\widehat\theta^2|F^\varepsilon_{1,x}|^2{\rm d}x{\rm d}y{\rm d}t\nonumber\\
&-\frac{1}{2}\mathbb{E}\int_{Q_T}(x+\varepsilon)^{2\gamma} l_y^2 \widehat\theta^2|F^\varepsilon_1|^2{\rm d}x{\rm d}y{\rm d}t
-2\sigma \mathbb{E}\int_{Q_T} l_x^2\widehat\theta^2 |F^\varepsilon_1|^2{\rm d}x{\rm d}y{\rm d}t\nonumber\\
&-2\sigma \mathbb{E}\int_{Q_T}\widehat\theta^2|F^\varepsilon_{1,x}|^2{\rm d}x{\rm d}y{\rm d}t\nonumber\\
\geq&(\frac{1}{2}-\epsilon-2\sigma) \mathbb{E}\int_{Q_T}\widehat\theta^2|F^\varepsilon_{1,x}|^2{\rm d}x{\rm d}y{\rm d}t-C(\epsilon) \mathbb{E}\int_{Q_T}s^2\lambda^2\widehat\phi^2\xi^2\widehat\theta^2|F^\varepsilon_1|^2{\rm d}x{\rm d}y{\rm d}t.
\end{align}
Taking $\epsilon=\frac{1}{2}-2\sigma>0$ due to $0\leq\sigma<\frac{1}{4}$, from (\ref{3.1.29}) we deduce (\ref{3.1.30}).\hfill$\Box$

\vspace{2mm}

Combining Lemma 3.3-Lemma 3.5, we have the following result.

\begin{lem} There exists $C=C(I,T,\gamma,\sigma,\mu,M,\lambda)$  such that
\begin{align}\label{1-3.35}
&\mathbb{E}\int_{Q_T} s\lambda^2\xi\widehat\theta^2 |v^\varepsilon_x|^2{\rm d}x{\rm
d}y{\rm d}t+\mathbb{E}\int_{Q_T}s\lambda^2(x+\varepsilon)^{2\gamma}\xi \widehat\theta^2|v^\varepsilon_y|^2{\rm d}x{\rm d}y{\rm d}t\nonumber\\
&+\mathbb{E}\int_{Q_T}s^3\lambda^4 \xi^3\widehat\theta^2|v^\varepsilon|^2{\rm d}x{\rm d}y{\rm d}t\nonumber\\
\leq &C\mathbb{E}\int_{Q_T}\widehat\theta^2|f_1|^2{\rm d}x{\rm d}y{\rm d}t+C\mathbb{E}\int_{Q_T}s^2\lambda^2\xi^2\widehat\theta^2|F^\varepsilon_1|^2{\rm d}x{\rm d}y{\rm d}t\nonumber\\
&+C\mathbb{E}\int_{\omega^{(2)}_T}\widehat\theta^2\left[s\xi|v^\varepsilon_x|^2+s(x+\varepsilon)^{2\gamma}\xi|v^\varepsilon_y|^2+s^3\xi^3|v^\varepsilon|^2\right]{\rm d}x{\rm d}y{\rm d}t
\end{align}
for all large $\lambda$ and $s$.
\end{lem}

{\noindent\bf Proof.}\ By substituting (\ref{3.15}), (\ref{1-3.27}) and (\ref{3.1.30}) into (\ref{3.1.18}), we find that
\begin{align}\label{1-3.36}
&\mathbb{E}\int_{Q_T} P_2\widehat\theta \left[{\rm d}v^\varepsilon+v^\varepsilon_{xx}{\rm d}t+\left(x+\varepsilon\right)^{2\gamma}v^\varepsilon_{yy}{\rm d}t+\frac{\sigma}{(x+\varepsilon)^2}v^\varepsilon{\rm d}t\right]{\rm d}x{\rm d}y\nonumber\\
\geq &\frac{1}{2}\mathbb{E}\int_{Q_T}|P_2|^2{\rm d}x{\rm d}y{\rm d}t-\frac{1}{2}\mathbb{E}\int_{Q_T}|P|^2{\rm d}x{\rm d}y{\rm d}t\nonumber\\
&+C\mathbb{E}\int_{Q_T}s^3\lambda^4\widehat\phi^3 \xi^3|z|^2{\rm d}x{\rm d}y{\rm d}t+C\mathbb{E}\int_{Q_T} s\lambda^2\widehat\phi \xi |z_x|^2{\rm d}x{\rm
d}y{\rm d}t\nonumber\\
&+C\mathbb{E}\int_{Q_T}s\lambda^2(x+\varepsilon)^{2\gamma}\widehat\phi \xi |z_y|^2{\rm d}x{\rm d}y{\rm d}t-C \mathbb E\int_{\Gamma_T}s\lambda \widehat \phi\xi |z_x|^2{\rm d}y{\rm d}t\nonumber\\
&-C\mathbb{E}\int_{Q_T}s^2\lambda^2\widehat\phi^2\xi^2\widehat\theta^2|F^\varepsilon_1|^2{\rm d}x{\rm d}y{\rm d}t
\end{align}
for all large $s\geq s_1$, $\lambda\geq \lambda_1$ and $z=\widehat \theta v^\varepsilon$.

In order to eliminate the boundary term, we introduce a cut-function $\chi\in C^2(\overline I)$
such that
\begin{align}\label{3.1.47}\left\{\begin{array}{ll}\chi(x,y)=0,&(x,y)\in \overline {\omega^{(1)}},\\
0<\chi(x,y)<1,&(x,y)\in \omega^{(2)}\backslash\omega^{(1)},\\
\chi(x,y)=1,&(x,y)\in \overline{I\backslash \omega^{(2)}}.\end{array}\right.
\end{align}
Then $\tilde v^\varepsilon:=\chi v^\varepsilon$ satisfies
\begin{align}\label{3.1.48}
\hspace{-0.1cm}\left\{\begin{array}{ll}
{\rm d}\tilde v^\varepsilon+\tilde v^\varepsilon_{xx}{\rm d}t+\left(x+\varepsilon\right)^{2\gamma}\tilde v^\varepsilon_{yy}{\rm d}t+\frac{\sigma}{(x+\varepsilon)^2}\tilde v^\varepsilon{\rm d}t=\tilde f_1{\rm d}t+\tilde F_1^\varepsilon{\rm d}B(t),&(x,y,t)\in Q_T,\\
\tilde v^\varepsilon(x,y,t)=0,& (x,y,t)\in \Sigma_T,\\
\tilde v^\varepsilon(x,y,T)=\chi(x,y) v^\varepsilon_T(x,y), &(x,y)\in I,
\end{array}\right.
\end{align}
where \begin{align*}
\tilde f_1=2\chi_x v_x^\varepsilon+\chi_{xx}v^\varepsilon+(x+\varepsilon)^{2\gamma}\left(2\chi_y v^\varepsilon_{y}+\chi_{yy}v^\varepsilon\right)+\chi f_1,\quad \tilde F_1^\varepsilon=\chi F_1^\varepsilon.
\end{align*}

Let $\tilde z=\widehat\theta \tilde v^\varepsilon$ and $\tilde P_2$, $\tilde P$ denote the same expressions as $P_2$, $P$ by replacing $z$ with $\tilde z$. By the definition of $\chi$, we obtain $\tilde z_x=0$ on $\Gamma_T$. Then applying (\ref{1-3.36}) to $\tilde v^\varepsilon$ yields
\begin{align}\label{}
&\mathbb{E}\int_{Q_T}\tilde P_2 \widehat\theta\left[{\rm d}\tilde  v^\varepsilon+\tilde v^\varepsilon_{xx}{\rm d}t+\left(x+\varepsilon\right)^{2\gamma}\tilde v^\varepsilon_{yy}{\rm d}t+\frac{\sigma}{(x+\varepsilon)^2}\tilde v^\varepsilon{\rm d}t\right]{\rm d}x{\rm d}y\nonumber\\
\geq &\frac{1}{2}\mathbb{E}\int_{Q_T}|\tilde P_2|^2{\rm d}x{\rm d}y{\rm d}t-\frac{1}{2}\mathbb{E}\int_{Q_T}|\tilde P|^2{\rm d}x{\rm d}y{\rm d}t\nonumber\\
&+C\mathbb{E}\int_{Q_T}s^3\lambda^4\widehat\phi^3 \xi^3|\tilde z|^2{\rm d}x{\rm d}y{\rm d}t+C\mathbb{E}\int_{Q_T} s\lambda^2\widehat\phi\xi |\tilde z_x|^2{\rm d}x{\rm
d}y{\rm d}t\nonumber\\
&+C\mathbb{E}\int_{Q_T}s\lambda^2(x+\varepsilon)^{2\gamma}\widehat\phi \xi |\tilde z_y|^2{\rm d}x{\rm d}y{\rm d}t-C\mathbb{E}\int_{Q_T}s^2\lambda^2\widehat\phi^2\xi^2\widehat\theta^2|\tilde F^\varepsilon_1|^2{\rm d}x{\rm d}y{\rm d}t
\end{align}
for all large $\lambda$ and $s$, which implies
\begin{align}\label{1-3.40}
&\mathbb{E}\int_{Q_T} s\lambda^2\widehat\phi\xi |\tilde z_x|^2{\rm d}x{\rm
d}y{\rm d}t+\mathbb{E}\int_{Q_T}s\lambda^2(x+\varepsilon)^{2\gamma}\widehat\phi \xi |\tilde z_y|^2{\rm d}x{\rm d}y{\rm d}t\nonumber\\
&+\mathbb{E}\int_{Q_T}s^3\lambda^4\widehat\phi^3 \xi^3|\tilde z|^2{\rm d}x{\rm d}y{\rm d}t+\mathbb{E}\int_{Q_T}|\tilde P_2|^2{\rm d}x{\rm d}y{\rm d}t\nonumber\\
\leq &C\mathbb{E}\int_{Q_T} \tilde P_2 \widehat\theta\left[{\rm d}\tilde  v^\varepsilon+\tilde v^\varepsilon_{xx}{\rm d}t+\left(x+\varepsilon\right)^{2\gamma}\tilde v^\varepsilon_{yy}{\rm d}t+\frac{\sigma}{(x+\varepsilon)^2}\tilde v^\varepsilon{\rm d}t\right]{\rm d}x{\rm d}y\nonumber\\
&+C\mathbb{E}\int_{Q_T}|\tilde P|^2{\rm d}x{\rm d}y{\rm d}t+C\mathbb{E}\int_{Q_T}s^2\lambda^2\widehat\phi^2\xi^2\widehat\theta^2| F^\varepsilon_1|^2{\rm d}x{\rm d}y{\rm d}t.
\end{align}
Using the equation of $\tilde v^\varepsilon$, ${\rm
Supp}(\chi_x), {\rm Supp}(\chi_y)\subset\omega^{(2)}$ and noticing that
\begin{align*}
\mathbb{E}\int_{Q_T} \tilde P_2 \widehat\theta \tilde F^\varepsilon_1{\rm d}x{\rm d}y{\rm d}B(t)=0,
\end{align*}
we see that
\begin{align}\label{3.1.43}
&\mathbb{E}\int_{Q_T}\tilde P_2\widehat\theta\left[{\rm d}\tilde  v^\varepsilon+\tilde v^\varepsilon_{xx}{\rm d}t+\left(x+\varepsilon\right)^{2\gamma}\tilde v^\varepsilon_{yy}{\rm d}t+\frac{\sigma}{(x+\varepsilon)^2}\tilde v^\varepsilon{\rm d}t\right]{\rm d}x{\rm d}y\nonumber\\
=&\mathbb{E}\int_{Q_T}\tilde P_2\widehat\theta\tilde f_1{\rm d}x{\rm d}y{\rm d}t+\mathbb{E}\int_{Q_T}\tilde P_2\widehat\theta\tilde F^\varepsilon_1{\rm d}x{\rm d}y{\rm d}B(t)\nonumber\\
\leq& \epsilon\mathbb{E}\int_{Q_T}|\tilde P_2|^2{\rm d}x{\rm d}y{\rm d}t+C(\epsilon)\mathbb{E}\int_{Q_T}\widehat\theta^2|f_1|^2{\rm d}x{\rm d}y{\rm d}t\nonumber\\
&+C(\epsilon)\mathbb{E}\int_{\omega^{(2)}_T}\widehat\theta^2\left[|v^\varepsilon_x|^2+(x+\varepsilon)^{2\gamma}|v^\varepsilon_y|^2+|v^\varepsilon|^2\right]{\rm d}x{\rm d}y{\rm d}t.
\end{align}
From \eqref{3.1.4}, \eqref{3.1.5} and \eqref{3.1.19}, we have
\begin{align}\label{3.1.45}
\mathbb{E}\int_{Q_T} |\tilde P|^2{\rm d}x{\rm d}y{\rm d}t=&\mathbb{E}\int_{Q_T}\left|(\tau-1) l_{xx}\tilde z- l_t\tilde z-(x+\varepsilon)^{2\gamma} l_{yy}\tilde z\right|^2{\rm d}x{\rm d}y{\rm d}t\nonumber\\
\leq &C(\lambda)\mathbb{E}\int_{Q_T}s^2\widehat\phi^2\xi^3 |\tilde z|^2{\rm d}x{\rm d}y{\rm d}t.
\end{align}
Substituting (\ref{3.1.43}) and (\ref{3.1.45}) into (\ref{1-3.40}) and choosing $\epsilon$ sufficiently small and $s$ sufficiently large, we then obtain
\begin{align}
&\mathbb{E}\int_{Q_T} s\lambda^2 \xi |\tilde z_x|^2{\rm d}x{\rm
d}y{\rm d}t+\mathbb{E}\int_{Q_T}s\lambda^2(x+\varepsilon)^{2\gamma}\widehat\phi\xi |\tilde z_y|^2{\rm d}x{\rm d}y{\rm d}t\nonumber\\
&+\mathbb{E}\int_{Q_T}s^3\lambda^4\widehat\phi^3 \xi^3|\tilde z|^2{\rm d}x{\rm d}y{\rm d}t\nonumber\\
\leq &C\mathbb{E}\int_{Q_T}\widehat\theta^2|f_1|^2{\rm d}x{\rm d}y{\rm d}t+C\mathbb{E}\int_{Q_T}s^2\lambda^2\widehat\phi^2\xi^2\widehat\theta^2| F^\varepsilon_1|^2{\rm d}x{\rm d}y{\rm d}t\nonumber\\
&+C\mathbb{E}\int_{\omega^{(2)}_T}\widehat\theta^2\left[|v^\varepsilon_x|^2+(x+\varepsilon)^{2\gamma}|v^\varepsilon_y|^2+|v^\varepsilon|^2\right]{\rm d}x{\rm d}y{\rm d}t.
\end{align}
 Using $\tilde z=z$ on $\overline{I\setminus\omega^{(2)}}$, we further have
\begin{align}
&\mathbb{E}\int_{Q_T\setminus \omega^{(2)}_T} s\lambda^2\widehat\phi\xi |z_x|^2{\rm d}x{\rm
d}y{\rm d}t+\mathbb{E}\int_{Q_T\setminus \omega^{(2)}_T}s\lambda^2(x+\varepsilon)^{2\gamma}\widehat\phi \xi |z_y|^2{\rm d}x{\rm d}y{\rm d}t\nonumber\\
&+\mathbb{E}\int_{Q_T\setminus \omega^{(2)}_T}s^3\lambda^4\widehat\phi^3 \xi^3|z|^2{\rm d}x{\rm d}y{\rm d}t\nonumber\\
\leq &C\mathbb{E}\int_{Q_T}\widehat\theta^2|f_1|^2{\rm d}x{\rm d}y{\rm d}t+C\mathbb{E}\int_{Q_T}s^2\lambda^2\widehat\phi^2\xi^2\widehat\theta^2| F^\varepsilon_1|^2{\rm d}x{\rm d}y{\rm d}t\nonumber\\
&+C(\lambda)\mathbb{E}\int_{\omega^{(2)}_T}\widehat\theta^2\left[|v^\varepsilon_x|^2+(x+\varepsilon)^{2\gamma}|v^\varepsilon_y|^2+|v^\varepsilon|^2\right]{\rm d}x{\rm d}y{\rm d}t.
\end{align}

Finally using $z=\widehat\theta v^\varepsilon$ and going back to $v^\varepsilon$, we obtain (\ref{1-3.35}) and complete the proof of this lemma. \hfill$\Box$

\vspace{2mm}

In order to prove Theorem 3.1, we also need the following the Cacciopoli
inequality.

\begin{lem}  Let $\gamma>0$, $0\leq\sigma<\frac{1}{4}$, $v_T\in L^2(\Omega,\mathcal F_T, \mathbb P; L^2(I))$, $f_1\in L^2_{\mathcal F}(0,T; L^2(I))$, $F^\varepsilon_1\in L^2_{\mathcal F}(0,T;L^2(I))$. Then there exists constant $C=C(I,T,\gamma,\omega,\sigma,\mu,M,\lambda)$  such that the solution $u^\varepsilon\in \mathcal H_T$ of the backward stochastic Grushin equation (\ref{3.1.2}) satisfies
\begin{align}\label{1-3.45}
&\mathbb{E}\int_{\omega^{(2)}_T}\widehat\theta^2\left[s\xi|v^\varepsilon_x|^2+s(x+\varepsilon)^{2\gamma}\xi|v^\varepsilon_y|^2\right]{\rm
d}x{\rm d}y{\rm d}t\nonumber\\
\leq &C\mathbb{E}\int_{\omega_T}s^3\xi^3\widehat\theta^2|v^\varepsilon|^2{\rm
d}x{\rm d}y{\rm d}t+\mathbb{E}\int_{Q_T}\widehat\theta^2|f_1|^2{\rm d}x{\rm d}y{\rm d}t.
\end{align}
\end{lem}

{\noindent\bf Proof.}\
 We choose a cut-function
$\zeta\in C^2(\overline I)$ such that $0\leq \zeta\leq 1$ and
$\zeta=1$ in $\omega^{(2)}$, $\zeta=0$ in $I\setminus \omega$. By It\^{o} formula, we have
\begin{align}
{\rm d}\big(\xi\widehat\theta^2|v^\varepsilon|^2\big)=\big(\xi_t\widehat\theta^2+2\xi\widehat\theta\widehat\theta_t\big)|v^\varepsilon|^2{\rm d}t+2\xi\widehat\theta^2v^\varepsilon{\rm d}v^\varepsilon+ \xi\widehat\theta^2({\rm d}v^\varepsilon)^2.
\end{align}
Together with the equation of $v^\varepsilon$ in (\ref{3.1.2}), we find that
\begin{align}\label{3.1.53}
0=&\mathbb{E}\int_{Q_T}\zeta^2{\rm d}\big(\xi\widehat\theta^2|v^\varepsilon|^2\big){\rm d}x{\rm d}y\nonumber\\
=&\mathbb{E}\int_{Q_T}\zeta^2\left[\big(\xi_t\widehat\theta^2+2\xi\widehat\theta\widehat\theta_t\big)|v^\varepsilon|^2{\rm d}t+2\xi\widehat\theta^2 v^\varepsilon{\rm d}v^\varepsilon+ \xi\widehat\theta^2({\rm d}v^\varepsilon)^2\right]{\rm d}x{\rm d}y\nonumber\\
=&\mathbb{E}\int_{Q_T}\zeta^2\xi\widehat\theta^2\bigg[\big(\xi^{-1}\xi_t+2s\widehat\varphi_t\big)|v^\varepsilon|^2
+2 v^\varepsilon\Big(-v^\varepsilon_{xx}-(x+\varepsilon)^{2\gamma}v^\varepsilon_{yy}\nonumber\\
&-\frac{\sigma}{(x+\varepsilon)^2}v^\varepsilon+f_1\Big)
+|F_1^\varepsilon|^2\bigg]{\rm d}x{\rm d}y{\rm d}t\nonumber\\
=&\mathbb{E}\int_{Q_T}\zeta^2\xi\widehat\theta^2\left[\big(\xi^{-1}\xi_t+2s\widehat\varphi_t\big)|v^\varepsilon|^2+2|v^\varepsilon_x|^2+2(x+\varepsilon)^{2\gamma}|v^\varepsilon_{y}|^2
-\frac{2\sigma}{(x+\varepsilon)^2}|v^\varepsilon|^2\right]{\rm d}x{\rm
d}y{\rm d}t\nonumber\\
&-\mathbb{E}\int_{Q_T}\xi\left(\zeta^2\widehat\theta^2\right)_{xx}|v^\varepsilon|^2{\rm d}x{\rm
d}y{\rm d}t-\mathbb{E}\int_{Q_T}(x+\varepsilon)^{2\gamma}\xi\left(\zeta^2
\widehat\theta^2\right)_{yy}|v^\varepsilon|^2{\rm
d}x{\rm d}y{\rm d}t\nonumber\\&+\mathbb{E}\int_{Q_T}\zeta^2\xi\widehat\theta^2\left(2f_1v^\varepsilon+|F^\varepsilon_1|^2\right){\rm d}x{\rm
d}y{\rm d}t,
\end{align}
which implies
\begin{align}\label{3.1.54}
&2\mathbb{E}\int_{Q_T}\zeta^2\xi\widehat\theta^2\left[|v^\varepsilon_x|^2+(x+\varepsilon)^{2\gamma}|v^\varepsilon_y|^2\right]{\rm d}x{\rm d}y{\rm d}t+\mathbb{E}\int_{Q_T}\zeta^2\xi\widehat\theta^2|F^\varepsilon_1|^2{\rm d}x{\rm d}y{\rm d}t\nonumber\\
\leq&\mathbb{E}\int_{Q_T}\xi\left[-\zeta^2\xi^{-1}\xi_t\widehat\theta^2-2s\zeta^2\widehat\varphi_t \widehat\theta^2 +\left(\zeta^2\widehat\theta^2\right)_{xx}+(x+\varepsilon)^{2\gamma}\left(\zeta^2
\widehat\theta^2\right)_{yy}+s\zeta^2\xi\widehat\theta^2\right]|v^\varepsilon|^2{\rm d}x{\rm
d}y{\rm d}t\nonumber\\
&+2\sigma\mathbb{E}\int_{Q_T}\frac{1}{(x+\varepsilon)^2}\zeta^2\xi\widehat\theta^2|v^\varepsilon|^2{\rm d}x{\rm d}y{\rm d}t+\mathbb{E}\int_{Q_T}s^{-1}\zeta^2\widehat\theta^2|f_1|^2{\rm d}x{\rm d}y{\rm d}t.
\end{align}
On the other hand, by Hardy inequality (\ref{1-1.5}), we have
\begin{align}\label{3.1.55}
&\sigma\mathbb{E}\int_{Q_T}\frac{1}{(x+\varepsilon)^2}\zeta^2\xi\widehat\theta^2|v^\varepsilon|^2{\rm
d}x{\rm d}y{\rm d}t\leq 4\sigma \mathbb{E}\int_{Q_T}\xi\left|\left(\zeta\widehat\theta v^\varepsilon\right)_x\right|^2{\rm
d}x{\rm d}y{\rm d}t\nonumber\\
\leq &4\sigma\mathbb{E}\int_{Q_T}\zeta^2\xi\widehat\theta^2|v^\varepsilon_x|^2{\rm d}x{\rm
d}y{\rm d}t+C(\lambda)\mathbb{E}\int_{Q_T}\left(\zeta^2+\zeta_x^2\right)s^2\xi^3\widehat\theta^2|v^\varepsilon|^2{\rm d}x{\rm d}y{\rm d}t.
\end{align}
Therefore, by the definition of $\zeta$ and
\begin{align*}
&\xi\left|-\zeta^2\xi^{-1}\xi_t\widehat\theta^2-2s\zeta^2\widehat\varphi_t \widehat\theta^2 +\left(\zeta^2\widehat\theta^2\right)_{xx}+(x+\varepsilon)^{2\gamma}\left(\zeta^2
\widehat\theta^2\right)_{yy}+s\zeta^2\xi\widehat\theta^2\right|\\
\leq & C(\lambda)s^2\left(\zeta^2+\zeta_x^2\right)\xi^3\widehat\theta^2,
\end{align*}
 we deduce from (\ref{3.1.54}) and (\ref{3.1.55}) that
\begin{align}\label{3.1.56}
&\mathbb{E}\int_{\omega^{(2)}_T}\xi\widehat\theta^2\left[(1-4\sigma)|v^\varepsilon_x|^2+(x+\varepsilon)^{2\gamma}|v^\varepsilon_y|^2\right]{\rm
d}x{\rm d}y{\rm d}t\nonumber\\
\leq &C(\lambda)\mathbb{E}\int_{\omega_T}s^2\xi^3\widehat\theta^2|v^\varepsilon|^2{\rm
d}x{\rm d}y{\rm d}t+\mathbb{E}\int_{Q_T}s^{-1}\widehat\theta^2|f_1|^2{\rm d}x{\rm d}y{\rm d}t.
\end{align}
Noticing that $0\leq\sigma<\frac{1}{4}$, by (\ref{3.1.56}) we immediately obtain (\ref{1-3.45}).\hfill$\Box$

\vspace{2mm}

Now we prove Theorem 3.1.

\vspace{2mm}

{\noindent\bf Proof of Theorem 3.1.}\ By Lemma 3.6 and Lemma 3.7, we have
\begin{align}\label{1-3.51}
&\mathbb{E}\int_{Q_T} s\lambda^2\xi\widehat\theta^2 |v^\varepsilon_x|^2{\rm d}x{\rm
d}y{\rm d}t+\mathbb{E}\int_{Q_T}s\lambda^2(x+\varepsilon)^{2\gamma}\xi \widehat\theta^2|v^\varepsilon_y|^2{\rm d}x{\rm d}y{\rm d}t\nonumber\\
&+\mathbb{E}\int_{Q_T}s^3\lambda^4 \xi^3\widehat\theta^2|v^\varepsilon|^2{\rm d}x{\rm d}y{\rm d}t\nonumber\\
\leq &C\mathbb{E}\int_{Q_T}\widehat\theta^2|f_1|^2{\rm d}x{\rm d}y{\rm d}t+C(\lambda)\mathbb{E}\int_{Q_T}s^2\xi^2\widehat\theta^2| F^\varepsilon_1|^2{\rm d}x{\rm d}y{\rm d}t\nonumber\\
&+C(\lambda)\mathbb{E}\int_{\omega_T}s^3\xi^3\widehat\theta^2|v^\varepsilon|^2{\rm d}x{\rm d}y{\rm d}t
\end{align}
for all large $\lambda$ and $s$. By a similar argument to (\ref{1-2.9}), we could prove $v^\varepsilon\rightarrow v$ in $\mathcal G_T$. Therefore, by letting $\varepsilon\rightarrow 0$ in
(\ref{1-3.51}), together with $F_1^\varepsilon\rightarrow F_1$ in $L^2_{\mathcal F}(0,T;L^2(I))$, we obtain (\ref{3.1.8}) and then complete the proof of Theorem 3.1.\hfill$\Box$

\subsection{Carleman estimate for forward stochastic Grushin equation with regular weight function}

In this subsection, we introduce a new regular weight function to establish the other Carlemen estimate for the forward stochastic Grushin equation with singular potential
\bes\label{3.2.57}\left\{\begin{array}{ll}
{\rm d}w-w_{xx}{\rm d}t-x^{2\gamma}w_{yy}{\rm d}t-\frac{\sigma}{x^2}w{\rm d}t=f_2{\rm d}t+F_2{\rm d}B(t),&(x,y,t)\in Q_T,\\
w(x,y,t)=0,& (x,y,t)\in \Sigma_T,\\
w(x,y,0)=0,&(x,y)\in I.\end{array}\right.
\ees
The regular weight function allows us to put the random source on the left-hand side of this Carleman estimate. Based on such a Carleman estimate we can obtain the uniqueness for our inverse problem.

We set
\begin{align}\label{3.2.59}
&\Phi(x,y,t)=e^{\lambda\varrho(x,y,t)},\quad \Theta (x,y,t)=e^{s\Phi(x,y,t)}
\end{align}
with
\begin{align}\label{3.2.60}
\varrho(x,y,t)=x^{2+2\gamma}y(1-y)-\mu x-(\lambda-t)^2+2\lambda^2.
\end{align}
Here $\mu$ is the same as the one in Section 3.1. We easily see that $\varrho>0$ in $Q_T$ if we choose $\lambda$ suitable large.

Our main result in this subsection is the following Carleman estimate for (\ref{3.2.57}) with regular weight function.

\begin{thm} Let $\gamma>0$, $0\leq\sigma<\frac{1}{4}$, $f_2\in L^2_{\mathcal F}(0,T; L^2(I))$, $F_2\in L^2_{\mathcal F}(0,T;H^1(I))$. Then there exist constants $\lambda_2=\lambda_2(I,T,\gamma,\sigma,\mu$), $s_2=s_2(I$, $T,\gamma,\sigma, \mu$, $\lambda)$,
$C=C(I,T,\gamma,\sigma,\mu)$ such that
\begin{align}\label{3.2.63}
&\mathbb{E}\int_{Q_T} s\lambda^2\Phi\Theta^2|w_x|^2{\rm d}x{\rm
d}y{\rm d}t+\mathbb{E}\int_{Q_T}s\lambda^2 \Phi \Theta^2x^{2\gamma} |w_y|^2{\rm d}x{\rm d}y{\rm d}t\nonumber\\
&+\mathbb{E}\int_{Q_T}s^3\lambda^4\Phi^3\Theta^2|w|^2{\rm d}x{\rm d}y{\rm d}t+\mathbb{E}\int_{Q_T}s\lambda\Phi\Theta^2|F_2|^2{\rm d}x{\rm d}y{\rm d}t\nonumber\\
\leq &C\mathbb{E}\int_{Q_T}\Theta^2 |f_2|^2{\rm d}x{\rm d}y{\rm d}t+C\mathbb{E}\int_{Q_T}s\Phi\Theta^2|\nabla F_{2}|^2{\rm d}x{\rm d}y{\rm d}t\nonumber\\
&+C\mathbb{E}\int_{I}s^2\lambda^2\Phi^2(T)\Theta^2(T)w^2(T){\rm d}x{\rm d}y+C\mathbb E\int_{\Gamma_T} s\lambda \Phi \Theta^2 |w_x|^2{\rm d}y{\rm d}t
\end{align}
for all $\lambda\geq\lambda_2$, $s\geq s_2$ and all $w\in \mathcal G_T$ satisfies \eqref{3.2.57}.
\end{thm}

\noindent{\bf Remark 3.1.}\  We could not eliminate the term of $ \nabla F_{2}$ on the right-hand side of (\ref{3.2.63}). Based on this reason, the random source $H$ to be determined in (\ref{1.3}) does not depend on spatial variables.

\vspace{2mm}

\noindent{\bf Remark 3.2.}\  The second large parameter $\lambda$ in studying the null controllability could be
omitted. However, it plays a very important role in determining the random source $H$. Therefore we have to separate $\lambda$ from constant $C$.

\vspace{2mm}

We still transfer to consider an approximate version of (\ref{3.2.57}):
\bes\label{3.2.58}\left\{\begin{array}{ll}
{\rm d}w^\varepsilon-w^\varepsilon_{xx}{\rm d}t-\left(x+\varepsilon\right)^{2\gamma}w^\varepsilon_{yy}{\rm d}t-\frac{\sigma}{(x+\varepsilon)^2}w^\varepsilon{\rm d}t=f_2{\rm d}t+F_2{\rm d}B(t),&(x,y,t)\in Q_T,\\
w^\varepsilon(x,y,t)=0,& (x,y,t)\in \Sigma_T,\\
w^\varepsilon(x,y,0)=0,&(x,y)\in I,\end{array}\right.
\ees
 where $0<\varepsilon<1$. Set
\begin{align*}
\widehat\Phi(x,y,t)=\Phi(x+\varepsilon,y,t),\quad \widehat\varrho(x,y,t)=\varrho(x+\varepsilon,y,t),\quad \widehat\Theta(x,y,t)=\Theta(x+\varepsilon,y,t).
\end{align*}

We first give a weighted identity
for the approximate problem (\ref{3.2.58}).

\begin{lem}Let $\tau$ be a constant such that $2<\tau<3$. Assume that $w^\varepsilon$ is an $H^2(\mathbb R^2)$-valued continuous semimartingale. Set $L=s\widehat\Phi$, $Z=\widehat \Theta w^\varepsilon$ and
\begin{align*}
{{Q}}_1=&{\rm d}Z+2 L_xZ_x{\rm d}t+2\left(x+\varepsilon\right)^{2\gamma} L_yZ_y{\rm d}t+\tau L_{xx}Z{\rm d}t,\\
{{Q}}_2=&-L_t Z-Z_{xx}-\left(x+\varepsilon\right)^{2\gamma}Z_{yy}-L_x^2 Z-\left(x+\varepsilon\right)^{2\gamma}L_y^2 Z-\frac{\sigma}{(x+\varepsilon)^2}Z,\nonumber\\
{{Q}}=&-(\tau-1) L_{xx}Z+\left(x+\varepsilon\right)^{2\gamma} L_{yy}Z.
\end{align*}
Then for a.e. $(x,y)\in \mathbb R^2$, it holds that
\begin{align}\label{1-3.57}
& Q_2 \widehat\Theta\left[{\rm d}w^\varepsilon-w^\varepsilon_{xx}{\rm d}t-\left(x+\varepsilon\right)^{2\gamma}w^\varepsilon_{yy}{\rm d}t-\frac{\sigma}{(x+\varepsilon)^2}w^\varepsilon{\rm d}t\right]\nonumber\\
=& |Q_2| ^2{\rm d}t+Q_2Q{\rm d}t+\sum_{i=1}^{5}\overline X_i{\rm d}t+{\rm d}\overline Y+\overline{\{\cdot\}}_x+\overline{\{\cdot\cdot\}}_y+\overline J, \quad \mathbb P-a.s.,
\end{align}
where
\begin{align*}
\overline X_1=&\left[(\tau+1) L_{xx}-(x+\varepsilon)^{2\gamma} L_{yy}\right]Z_x^2, \nonumber\\
\overline X_2=&\left[-2\gamma(x+\varepsilon)^{2\gamma-1}L_{x}+(\tau-1)(x+\varepsilon)^{2\gamma}L_{xx}+(x+\varepsilon)^{4\gamma} L_{yy}\right]Z_y^2,\nonumber\\
\overline X_3=&4\left[\gamma(x+\varepsilon)^{2\gamma-1} L_y+(x+\varepsilon)^{2\gamma} L_{xy}\right]Z_x Z_y,\nonumber\\
\overline X_4=&\left[(3-\tau) L^2_{x} L_{xx}+2\gamma(x+\varepsilon)^{2\gamma-1} L_x L_y^2+3(x+\varepsilon)^{4\gamma} L_y^2 L_{yy}\right]Z^2\nonumber\\
&+(x+\varepsilon)^{2\gamma}\left[4 L_xL_y L_{xy}+L_{x}^2 L_{yy}+
(1-\tau) L_{xx} L_y^2\right]Z^2 \nonumber\\
&+\left[(1-\tau)\frac{\sigma}{(x+\varepsilon)^2} L_{xx}-\frac{2\sigma}{(x+\varepsilon)^3} L_x
+\frac{\sigma}{(x+\varepsilon)^{2-2\gamma}} L_{yy}\right]Z^2,\nonumber\\
\overline X_5=&\left[\frac{1}{2} L_{tt}+(1-\tau)L_{xx} L_t+2 L_x L_{xt}+2(x+\varepsilon)^{2\gamma} L_y L_{yt}+(x+\varepsilon)^{2\gamma} L_{yy} L_t
\right.\nonumber\\&-\left.\frac{\tau}{2} L_{xxxx}-\frac{\tau}{2}(x+\varepsilon)^{2\gamma} L_{xxyy}\right]Z^2,\nonumber\\
\overline Y=&\frac{1}{2}Z_x^2+\frac{1}{2}(x+\varepsilon)^{2\gamma}Z_y^2-\frac{1}{2}\left[ L_t+ L_x^2+(x+\varepsilon)^{2\gamma} L_y^2+\frac{\sigma}{(x+\varepsilon)^2}\right]Z^2,\nonumber\\
\overline{\{\cdot\}}=&-Z_x {\rm d}Z+\bigg[-L_x L_tZ^2- L_x Z_x^2+(x+\varepsilon)^{2\gamma} L_x Z^2_y- L_x^3 Z^2-(x+\varepsilon)^{2\gamma} L_xL_y^2Z^2\nonumber\\
&-\frac{\sigma}{(x+\varepsilon)^2} L_x Z^2-2(x+\varepsilon)^{2\gamma}L_y Z_x Z_y-\tau  L_{xx}ZZ_x+\frac{\tau}{2} L_{xxx}Z^2 \bigg]{\rm d}t,\nonumber\\
\overline{\{\cdot\cdot\}}=&-(x+\varepsilon)^{2\gamma}Z_y{\rm d}Z+\left[-(x+\varepsilon)^{2\gamma} L_y L_t Z^2-2(x+\varepsilon)^{2\gamma} L_x Z_x Z_y+(x+\varepsilon)^{2\gamma} L_{y}Z_x^2\right.\nonumber\\
&\left.-(x+\varepsilon)^{4\gamma} L_{y}Z_y^2-(x+\varepsilon)^{2\gamma} L_{x}^2 L_y Z^2-(x+\varepsilon)^{4\gamma} L_y^3 Z^2-\frac{\sigma}{(x+\varepsilon)^{2-2\gamma}}L_y Z^2\right.\nonumber\\
&-\left.\tau(x+\varepsilon)^{2\gamma} L_{xx} Z Z_{y}+\frac{\tau}{2}(x+\varepsilon)^{2\gamma} L_{xxy}Z^2\right]{\rm d}t,\\
\overline J=&-\frac{1}{2}({\rm d}Z_x)^2-\frac{1}{2}(x+\varepsilon)^{2\gamma}({\rm d}Z_y)^2+\frac{1}{2}\left[ L_t+ L_x^2+(x+\varepsilon)^{2\gamma}L_y^2+\frac{\sigma}{(x+\varepsilon)^2}\right]({\rm d}Z)^2.
\end{align*}
\end{lem}
{\noindent\bf Proof.}\  Notice that $\widehat\Theta=e^{L}$, $L=s\widehat\Phi$ and $Z=\widehat\Theta w^{\varepsilon}$. Then we have
\begin{align*}
\widehat\Theta\left[{\rm d}w^{\varepsilon}-w_{xx}^{\varepsilon}{\rm d}t-\left(x+\varepsilon\right)^{2\gamma}w_{yy}^{\varepsilon}{\rm d}t-\frac{\sigma}{(x+\varepsilon)^2}w^{\varepsilon}{\rm d}t\right]=Q_1+(Q_2+Q){\rm d}t.
\end{align*}
Hence
\begin{align}\label{3.2.65}
\hspace{-0.05cm} Q_2\widehat\Theta\left[{\rm d}w^{\varepsilon}-w_{xx}^{\varepsilon}{\rm d}t-\left(x+\varepsilon\right)^{2\gamma}w_{yy}^{\varepsilon}{\rm d}t-\frac{\sigma}{(x+\varepsilon)^2}w^{\varepsilon}{\rm d}t\right]=Q_1 Q_2+|Q_2|^2{\rm d}t+Q_2 Q{\rm d}t.
\end{align}

We only need to deal with $-L_t  Z {Q}_1  $ in $Q_1Q_2$. The calculations of the other terms are similar to the ones in $P_1P_2$. Therefore, by using a similar argument similar to Lemma 3.2, together with
\begin{align*}
- L_t Z{Q}_1=&-L_t Z\left[{\rm d}Z+2 L_xZ_x{\rm d}t+2\left(x+\varepsilon\right)^{2\gamma} L_yZ_y{\rm d}t+\tau L_{xx}Z{\rm d}t\right]\nonumber\\
=&-\frac{1}{2}{\rm d}( L_t Z^2)+\frac{1}{2} L_{tt} Z^2{\rm d}t+\frac{1}{2} L_t({\rm d} Z)^2-( L_x  L_t Z^2)_x{\rm d}t+ (L_{xx} L_t+L_x L_{xt}) Z^2{\rm d}t\nonumber\\
&-\left[\left(x+\varepsilon\right)^{2\gamma} L_y L_t Z^2\right]_y{\rm d}t+\left(x+\varepsilon\right)^{2\gamma}(L_{yy} L_t+ L_y L_{yt})Z^2{\rm d}t-\tau L_{xx}L_{t}Z^2{\rm d}t,
\end{align*}
 we obtain (\ref{1-3.57}).\hfill$\Box$

\vspace{2mm}

Now, integrating both sides of (\ref{1-3.57}) in $Q_T$, taking mathematical expectation in $\Omega$, we obtain
\begin{align}\label{3.2.68}
&\mathbb{E}\int_{Q_T}  Q_2\widehat\Theta\left[{\rm d}w^{\varepsilon}-w_{xx}^{\varepsilon}{\rm d}t-\left(x+\varepsilon\right)^{2\gamma}w_{yy}^{\varepsilon}{\rm d}t-\frac{\sigma}{(x+\varepsilon)^2}w^{\varepsilon}{\rm d}t\right]{\rm d}x{\rm d}y\nonumber\\
\geq &\frac{1}{2}\mathbb{E}\int_{Q_T}|Q_2|^2{\rm d}x{\rm d}y{\rm d}t-\frac{1}{2}\mathbb{E}\int_{Q_T}|Q|^2{\rm d}x{\rm d}y{\rm d}t+\sum_{i=1}^{5}\mathbb{E}\int_{Q_T}\overline X_i{\rm d}x{\rm d}y{\rm d}t\nonumber\\
&+\mathbb{E}\int_{Q_T}{\rm d}\overline Y{\rm d}x{\rm d}y+\mathbb{E}\int_{Q_T}\big(\overline{\{\cdot\}}_x+\overline{\{\cdot\cdot\}}_y\big){\rm d}x{\rm d}y +\mathbb{E}\int_{Q_T}\overline J{\rm d}x{\rm d}y.
\end{align}

In the following we estimate the last four terms on the right-hand side of (\ref{3.2.68}).

\begin{lem} There exists constant $C=C(I,T,\gamma,\sigma,\mu)$  such that
\begin{align}\label{1-3.60}
\sum_{i=1}^{5}\mathbb{E}\int_{Q_T}\overline X_i{\rm d}x{\rm d}y{\rm d}t\geq &C\mathbb{E}\int_{Q_T}s^3\lambda^4\widehat\Phi^3|Z|^2{\rm d}x{\rm d}y{\rm d}t+C\mathbb{E}\int_{Q_T} s\lambda^2\widehat\Phi|Z_x|^2{\rm d}x{\rm
d}y{\rm d}t\nonumber\\
&+C\mathbb{E}\int_{Q_T}s\lambda^2(x+\varepsilon)^{2\gamma}\widehat\Phi |Z_y|^2{\rm d}x{\rm d}y{\rm d}t
\end{align}
for all large $\lambda$ and $s$.
\end{lem}

{\noindent\bf Proof.}\ For regular weight function $\widehat\varrho(x,y,t)=(x+\varepsilon)^{2+2\gamma}y(1-y)-\mu(x+\varepsilon)-(\lambda-t)^2+2\lambda^2$, we have the following properties of $\widehat\varrho$:
\begin{align}\label{3.2.69}\left\{\begin{array}{l}
\widehat\varrho_t=2(\lambda-t),\ \ \widehat\varrho_{tt}=-2,\ \ \widehat\varrho_{xt}=\widehat\varrho_{yt}=0,\\
\widehat\varrho_x<-\delta_0,\ \  \widehat\varrho_x\widehat\varrho_{xxx}\leq 0,\ \ |\widehat\varrho_{xxxx}|\leq C(x+\varepsilon)^{-2},\ \ |\widehat\varrho_{y}|+|\widehat\varrho_{yy}|\leq C(x+\varepsilon)^{2+2\gamma},\\
|\widehat\varrho_{xy}|+|\widehat\varrho_{xx}|+|\widehat\varrho_{yy}|+|\widehat\varrho_{xxy}|+|\widehat\varrho_{xyy}|+|\widehat\varrho_{xxyy}|\leq C(x+\varepsilon)^{2\gamma}.
\end{array}\right.\end{align}
Then, by a similar process to Lemma 3.3, we could obtain (\ref{1-3.60}) for all large $\lambda$ and $s$. \hfill$\Box$

\begin{lem} There exists constant $C=C(I,T,\gamma,\sigma,\mu)$ such that
\begin{align}\label{1-3.62}
&\mathbb{E}\int_{Q_T}{\rm d}\overline Y{\rm d}x{\rm d}y\geq-C\mathbb{E}\int_{I}s^2\lambda^2\widehat\Phi^2(T)Z^2(T){\rm d}x{\rm d}y.
\end{align}
\end{lem}
{\noindent\bf Proof.}\ By $Z|_{t=0}=0$, ${\mathbb P}$-a.s. in $I$, we have
\begin{align}\label{3.2.77}
&\mathbb{E}\int_{Q_T}{\rm d}\overline Y{\rm d}x{\rm d}y\nonumber\\
=&\frac{1}{2}\mathbb{E}\int_I\big[|Z_{x}(T)|^2+(x+\varepsilon)^{2\gamma}|Z_{y}(T)|^2\big]{\rm d}x{\rm d}y-\frac{1}{2}\mathbb{E}\int_I\bigg[L_t(T)+
L_{x}^2(T)+(x+\varepsilon)^{2\gamma} L_{y}^2(T)\nonumber\\
&+\frac{\sigma}{(x+\varepsilon)^2}\bigg]|Z(T)|^2{\rm d}x{\rm d}y\nonumber\\
\geq&\frac{1}{2}\mathbb{E}\int_{I}\left[|Z_{x}(T)|^2-\frac{\sigma}{(x+\varepsilon)^2}|Z(T)|^2\right]{\rm d}x{\rm d}y
-C\mathbb{E}\int_{I}s^2\lambda^2\widehat\Phi ^2(T)|Z(T)|^2{\rm d}x{\rm d}y,
\end{align}
together with $0\leq\sigma<\frac{1}{4}$, which implies (\ref{1-3.62}).
\hfill$\Box$

\begin{lem}  There exists constant $C=C(I,T,\gamma,\sigma,\mu)$  such that
\begin{align}\label{1-3.64}
&\mathbb{E}\int_{Q_T}\big(\overline{\{\cdot\}}_x+\overline{\{\cdot\cdot\}}_y\big){\rm d}x{\rm d}y\geq -C\mathbb E\int_{\Gamma_T} s\lambda \widehat \Phi |Z_x|^2{\rm d}y{\rm d}t.
\end{align}
\end{lem}

{\noindent\bf Proof.}\ Since $\widehat \Phi$ has the same property (\ref{3.1.39}) as $\widehat\varphi$ on the boundary of $I$. Therefore we immediately obtain the estimate (\ref{1-3.64}) for boundary term on the right-hand side of (\ref{3.2.68}). \hfill$\Box$

\begin{lem} There exists constant $C=C(I,T,\gamma,\sigma,\mu)$  such that
\begin{align}\label{1-3.67}
&\mathbb{E}\int_{Q_T}\overline J{\rm d}x{\rm d}y\geq C\mathbb{E}\int_{Q_T}s\lambda\widehat\Phi\widehat\Theta^2|F_2|^2{\rm d}x{\rm d}y{\rm d}t-C\mathbb{E}\int_{Q_T}s\widehat\Phi\widehat\Theta^2|\nabla F_{2}|^2{\rm d}x{\rm d}y{\rm d}t
\end{align}
for all large $\lambda$ and $s$.
\end{lem}
{\noindent\bf Proof.}\ It is easily to see that
\begin{align*}\left\{\begin{array}{l}
({\rm d}Z)^2=\widehat\Theta^2|F_2|^2{\rm d}t,\\
({\rm d}Z_x)^2= L_{x}^2\widehat\Theta^2|F_2|^2{\rm d}t+2 L_x\widehat\Theta^2F_2F_{2,x}{\rm d}t+\widehat\Theta^2|F_{2,x}|^2{\rm d}t,\\
({\rm d}Z_y)^2=L_y^2\widehat\Theta^2|F_2|^2{\rm d}t+2 L_y\widehat\Theta^2F_2F_{2,y}{\rm d}t+\widehat\Theta^2|F_{2,y}|^2{\rm d}t.\end{array}\right.
\end{align*}
Therefore, we have
\begin{align}
&\mathbb{E}\int_{Q_T}\overline J{\rm d}x{\rm d}y\nonumber\\
=&\frac{1}{2}\mathbb{E}\int_{Q_T}\left[L_t+\frac{\sigma}{(x+\varepsilon)^2}\right]\widehat\Theta^2 |F_2|^2{\rm d}x{\rm d}y{\rm d}t-\frac{1}{2}\mathbb{E}\int_{Q_T}\widehat\Theta^2\big[|F_{2,x}|^2+(x+\varepsilon)^{2\gamma}|F_{2,y}|^2\big]{\rm d}x{\rm d}y{\rm d}t\nonumber\\
&-\mathbb{E}\int_{Q_T}\big[L_x\widehat\Theta^2 F_2F_{2,x}+(x+\varepsilon)^{2\gamma}L_y\widehat \Theta^2 F_2F_{2,y}\big]{\rm d}x{\rm d}y{\rm d}t.
\end{align}
By $\widehat\varrho_t=2(\lambda-t)$ and $0\leq\sigma<\frac{1}{4}$, we have
\begin{align}\label{1-3.69}
\frac{1}{2}\mathbb{E}\int_{Q_T}\left[L_t+\frac{\sigma}{(x+\varepsilon)^2}\right]\widehat\Theta^2 |F_2|^2{\rm d}x{\rm d}y{\rm d}t\geq \mathbb{E}\int_{Q_T} s\lambda(\lambda-T)\widehat\Phi\widehat\Theta^2|F_2|^2{\rm d}x{\rm d}y{\rm d}t.
\end{align}
On the other hand, by Young's inequality with $\epsilon>0$, we obtain
\begin{align}\label{1-3.70}
&-\mathbb{E}\int_{Q_T}\left[L_x\widehat\Theta^2 F_2F_{2,x}+(x+\varepsilon)^{2\gamma}L_y\widehat \Theta^2 F_2F_{2,y}\right]{\rm d}x{\rm d}y{\rm d}t\nonumber\\
=& -\mathbb{E}\int_{Q_T}s\lambda\left[\widehat \varrho_x F_2 F_{2,x}+(x+\varepsilon)^{2\gamma}\widehat \varrho_y F_2 F_{2,y}\right]\widehat\Phi\widehat\Theta^2{\rm d}x{\rm d}y{\rm d}t\nonumber\\
\geq& -\epsilon \mathbb{E}\int_{Q_T}s\lambda^2\widehat\Phi\widehat\Theta^2|F_2|^2{\rm d}x{\rm d}y{\rm d}t-C(\epsilon)\mathbb{E}\int_{Q_T}s\widehat\Phi\widehat\Theta^2\left(|F_{2,x}|^2+|F_{2,y}|^2\right){\rm d}x{\rm d}y{\rm d}t.
\end{align}
Therefore, from (\ref{1-3.69}) and (\ref{1-3.70}) we deduce that
\begin{align}
\mathbb{E}\int_{Q_T}\overline J{\rm d}x{\rm d}y\geq&\mathbb{E}\int_{Q_T}s\lambda\left[(1-\epsilon)\lambda-T\right]\widehat\Phi\widehat\Theta^2|F_2|^2{\rm d}x{\rm d}y{\rm d}t\nonumber\\
&-C\mathbb{E}\int_{Q_T}s\widehat\Phi\left(|F_{2,x}|^2+|F_{2,y}|^2\right)\widehat\Theta^2{\rm d}x{\rm d}y{\rm d}t.
\end{align}
Then taking $\epsilon$ sufficiently small and $\lambda$ sufficiently large, we obtain (\ref{1-3.67}). \hfill$\Box$

\vspace{2mm}

Now we proof Theorem 3.8.

\vspace{2mm}

{\noindent\bf Proof of Theorem 3.8.}\ Substituting \eqref{1-3.60}, \eqref{1-3.62}, \eqref{1-3.64} and \eqref{1-3.67} into \eqref{3.2.68}, we find that
\begin{align}\label{3.2.82}
&\mathbb{E}\int_{Q_T}  Q_2\widehat\Theta\left[{\rm d}w^{\varepsilon}-w_{xx}^{\varepsilon}{\rm d}t-\left(x+\varepsilon\right)^{2\gamma}w_{yy}^{\varepsilon}{\rm d}t-\frac{\sigma}{(x+\varepsilon)^2}w^{\varepsilon}{\rm d}t\right]{\rm d}x{\rm d}y\nonumber\\
\geq &\frac{1}{2}\mathbb{E}\int_{Q_T}|Q_2|^{2}{\rm d}x{\rm d}y{\rm d}t-\frac{1}{2}\mathbb{E}\int_{Q_T}|Q|^2{\rm d}x{\rm d}y{\rm d}t\nonumber\\
&+C\mathbb{E}\int_{Q_T}s^3\lambda^4\widehat\Phi^3|Z|^2{\rm d}x{\rm d}y{\rm d}t+C\mathbb{E}\int_{Q_T} s\lambda^2\widehat\Phi|Z_x|^2{\rm d}x{\rm
d}y{\rm d}t\nonumber\\
&+C\mathbb{E}\int_{Q_T}s\lambda^2(x+\varepsilon)^{2\gamma}\widehat\Phi |Z_y|^2{\rm d}x{\rm d}y{\rm d}t+C\mathbb{E}\int_{Q_T}s\lambda\widehat\Phi\widehat\Theta^2|F_2|^2{\rm d}x{\rm d}y{\rm d}t\nonumber\\
&-C\mathbb{E}\int_{Q_T}s\widehat\Phi\widehat\Theta^2|\nabla F_{2}|^2{\rm d}x{\rm d}y{\rm d}t
-C\mathbb{E}\int_{I}s^2\lambda^2\widehat\Phi^2(T)Z^2(T){\rm d}x{\rm d}y\nonumber\\
&-C\mathbb E\int_{\Gamma_T} s\lambda \widehat \Phi |Z_x|^2{\rm d}y{\rm d}t,
\end{align}
which implies
\begin{align}\label{1-3.73}
&\mathbb{E}\int_{Q_T} s\lambda^2\widehat\Phi|Z_x|^2{\rm d}x{\rm
d}y{\rm d}t+\mathbb{E}\int_{Q_T}s\lambda^2(x+\varepsilon)^{2\gamma}\widehat\Phi |Z_y|^2{\rm d}x{\rm d}y{\rm d}t\nonumber\\
&+\mathbb{E}\int_{Q_T}s^3\lambda^4\widehat\Phi^3|Z|^2{\rm d}x{\rm d}y{\rm d}t+\mathbb{E}\int_{Q_T}s\lambda\widehat\Phi\widehat\Theta^2|F_2|^2{\rm d}x{\rm d}y{\rm d}t+\frac{1}{2}\mathbb{E}\int_{Q_T}|Q_2|^{2}{\rm d}x{\rm d}y{\rm d}t\nonumber\\
\leq &C\mathbb{E}\int_{Q_T}  Q_2\widehat\Theta\left[{\rm d}w^{\varepsilon}-w_{xx}^{\varepsilon}{\rm d}t-\left(x+\varepsilon\right)^{2\gamma}w_{yy}^{\varepsilon}{\rm d}t-\frac{\sigma}{(x+\varepsilon)^2}w^{\varepsilon}{\rm d}t\right]{\rm d}x{\rm d}y\nonumber\\
&+C\mathbb{E}\int_{Q_T}|Q|^2{\rm d}x{\rm d}y{\rm d}t+C\mathbb{E}\int_{Q_T}s\widehat\Phi\widehat\Theta^2|\nabla F_{2}|^2{\rm d}x{\rm d}y{\rm d}t
+C\mathbb{E}\int_{I}s^2\lambda^2\widehat\Phi^2(T)Z^2(T){\rm d}x{\rm d}y\nonumber\\
&+C\mathbb E\int_{\Gamma_T} s\lambda \widehat \Phi |Z_x|^2{\rm d}y{\rm d}t.
\end{align}
Using the equation of $w^\varepsilon$ and noting that
$\mathbb{E}\int_{Q_T}  Q_2\widehat\Theta F_2{\rm d}x{\rm d}y{\rm d}B(t)=0$,
we have
\begin{align}\label{3.2.83}
&\mathbb{E}\int_{Q_T}  Q_2\widehat\Theta \left[{\rm d}w^{\varepsilon}-w_{xx}^{\varepsilon}{\rm d}t-\left(x+\varepsilon\right)^{2\gamma}w_{yy}^{\varepsilon}{\rm d}t-\frac{\sigma}{(x+\varepsilon)^2}w^{\varepsilon}{\rm d}t\right]{\rm d}x{\rm d}y\nonumber\\
\leq& \frac{1}{2}\mathbb{E}\int_{Q_T}|Q_2|^2{\rm d}x{\rm d}y{\rm d}t+\frac{1}{2}\mathbb{E}\int_{Q_T}\widehat\Theta^2 |f_2|^2{\rm d}x{\rm d}y{\rm d}t.
\end{align}
Obviously,
\begin{align}\label{3.2.80}
& |{Q}|^2\leq Cs^2\lambda^4\widehat\Phi^2|Z|^2.
\end{align}
From \eqref{1-3.73}-\eqref{3.2.80}, it follows that
\begin{align}\label{1-3.76}
&\mathbb{E}\int_{Q_T} s\lambda^2\widehat\Phi|Z_x|^2{\rm d}x{\rm
d}y{\rm d}t+\mathbb{E}\int_{Q_T}s\lambda^2(x+\varepsilon)^{2\gamma}\widehat\Phi |Z_y|^2{\rm d}x{\rm d}y{\rm d}t\nonumber\\
&+\mathbb{E}\int_{Q_T}s^3\lambda^4\widehat\Phi^3|Z|^2{\rm d}x{\rm d}y{\rm d}t+\mathbb{E}\int_{Q_T}s\lambda\widehat\Phi\widehat\Theta^2|F_2|^2{\rm d}x{\rm d}y{\rm d}t\nonumber\\
\leq &C\mathbb{E}\int_{Q_T}\widehat\Theta^2 |f_2|^2{\rm d}x{\rm d}y{\rm d}t+C\mathbb{E}\int_{Q_T}s\widehat\Phi\widehat\Theta^2|\nabla F_{2}|^2{\rm d}x{\rm d}y{\rm d}t\nonumber\\
&+C\mathbb{E}\int_{I}s^2\lambda^2\widehat\Phi^2(T)Z^2(T){\rm d}x{\rm d}y+C\mathbb E\int_{\Gamma_T} s\lambda \widehat \Phi |Z_x|^2{\rm d}y{\rm d}t
\end{align}
for all large $\lambda$ and $s$.

Finally, going back to the original variable $w^{\varepsilon}$ and letting $\varepsilon\rightarrow 0$ in (\ref{1-3.76}), we can obtain the desired estimate (\ref{3.2.63}). This completes the proof of Theorem 3.8. \hfill$\Box$

\section{Proof of Theorem 1.1}
\setcounter{equation}{0}
This section is devoted to proving the null controllability result for the forward stochastic Grushin equation \eqref{1.2}, i.e. Theorem 1.1.

\vspace{2mm}

{\noindent\bf Proof.}\  It is well known that the key ingredient for proving Theorem 1.1 is to obtain the observability
inequality for the corresponding adjoint system
\bes\label{1.4}\left\{\begin{array}{ll}
{\rm d}v+v_{xx}{\rm d}t+x^{2\gamma}v_{yy}{\rm d}t+\frac{\sigma}{x^2}v{\rm d}t=(-\alpha v-\beta V){\rm d}t+V{\rm d}B(t),&(x,y,t)\in Q_T,\\
v(x,y,t)=0,&(x,y,t)\in \Sigma_T,\\
v(x,y,T)=v_T(x,y),&(x,y)\in I.
\end{array}\right.
\ees
 More precisely, we will prove the following observability inequality for (\ref{1.4}):
\begin{align}\label{1.5}
&\mathbb E\int_{I}|v(0)|^2{\rm d}x{\rm d}y\leq C\mathbb E\int_{\omega_T}|v|^2{\rm d}x{\rm d}y{\rm d}t+C\mathbb E\int_{Q_T}|V|^2{\rm d}x{\rm d}y{\rm d}t,
\end{align}
where $C$ is depending on $I,T,\gamma,\omega,\sigma,\alpha$ and $\beta$.

  We apply Theorem 3.1 to \eqref{1.4} to obtain
\begin{align}\label{4.1}
&\mathbb{E}\int_{Q_T} s^3\xi^3\theta^2|v|^2{\rm d}x{\rm d}y{\rm d}t\leq C\mathbb{E}\int_{{\omega}_T}s^3\xi^3\theta^2|v|^2{\rm d}x{\rm d}y{\rm d}t+C\mathbb{E}\int_{Q_T}s^2\xi^2 \theta^2|V|^2{\rm d}x{\rm d}y{\rm d}t
\end{align}
for all large $\lambda>\lambda_1$ and $s>s_1$. We fix $\lambda=\lambda_1$ and $s=s_1$. By
$$M_1:=\max_{Q_T}(\xi^2+\xi^3)\theta^2<+\infty,$$
we further obtain
\begin{align}\label{4.3}
\mathbb{E}\int_{Q_T}\xi^3\theta^2 |v|^2{\rm d}x{\rm d}y{\rm
d}t\leq &C\mathbb{E}\int_{Q_T}\xi^2\theta^2|V|^2{\rm d}x{\rm d}y{\rm d}t
+\mathbb{E}\int_{\omega_T}\xi^3\theta^2|v|^2{\rm d}x{\rm d}y{\rm d}t\nonumber\\
\leq &C(\lambda_1,s_1,M_1)\left(\mathbb{E}\int_{Q_T}|V|^2{\rm d}x{\rm d}y{\rm d}t+\mathbb{E}\int_{{\omega_T}}|v|^2{\rm d}x{\rm d}y{\rm d}t\right).
\end{align}
On the other hand, by
$$m_1:=\min_{I\times (\frac{T}{4},\frac{3T}{4})}\xi^3\theta^2>0,$$
we obtain
\begin{align}\label{4.4}
&\mathbb{E}\int_{Q_T}\xi^3\theta^2|v|^2{\rm d}x{\rm d}y{\rm d}t\geq \mathbb{E}\int_{\frac{T}{4}}^{\frac{3T}{4}}\int_{I}\xi^3\theta^2|v|^2{\rm d}x{\rm d}y{\rm d}t\geq m_1\mathbb{E}\int_{\frac{T}{4}}^{\frac{3T}{4}}\int_{I}|v|^2{\rm d}x{\rm d}y{\rm d}t.
\end{align}
From \eqref{4.3} and \eqref{4.4}, we deduce
\begin{align}\label{4.5}
&\mathbb{E}\int_{\frac{T}{4}}^{\frac{3T}{4}}\int_{I}|v|^2{\rm d}x{\rm d}y{\rm d}t\leq  C(\lambda_1,s_1,M_1,m_1)\left(\mathbb{E}\int_{{\omega_T}}|v|^2{\rm d}x{\rm d}y{\rm d}t+\mathbb{E}\int_{Q_T}|V|^2{\rm d}x{\rm d}y{\rm d}t\right).
\end{align}

By the standard estimate for the backward stochastic equation \eqref{1.4}, we obtain for any $0\leq \tau<\tilde \tau\leq T$ that
\begin{align}\label{4.6}
&\mathbb{E}\int_{I} |v(\tau)|^2{\rm d}x{\rm d}y\leq \mathbb{E}\int_{I} |v(\tilde\tau)|^2{\rm d}x{\rm d}y+C\mathbb{E}\int_{\tau}^{\tilde\tau}\int_{I}|v|^2{\rm d}x{\rm d}y{\rm d}t.
\end{align}
Then from the Gronwall inequality, it follows that
\begin{align}\label{4.7}
\mathbb{E}\int_{I} |v(\tau)|^2{\rm d}x{\rm d}y\leq e^{CT}\mathbb{E}\int_{I} |v(\tilde\tau)|^2{\rm d}x{\rm d}y,\quad 0\leq\tau<\tilde\tau\leq T.
\end{align}
Further, letting $\tau=0$ and integrating \eqref{4.7} over $(\frac{T}{4},\frac{3T}{4})$ with respect to $\tilde \tau$, we obtain
\begin{align}\label{4.8}
&\mathbb{E}\int_I |v(0)|^2{\rm d}x{\rm d}y\leq C\mathbb{E}\int_{\frac{T}{4}}^{\frac{3T}{4}}\int_I |v|^2{\rm d}x{\rm d}y{\rm d}t.
\end{align}
Combining \eqref{4.5} and \eqref{4.8}, we obtain \eqref{1.5}. Then by a standard dual argument, e.g. as [\ref{Zhang2009}] or [\ref{Yan2018JMAA}], we could obtain a pair $(g, G)\in L_\mathcal F^2(0, T; L^2(\omega))\times L_\mathcal F^2(0, T; L^2(I))$ that drives the
corresponding solution $u$ of (\ref{1.2}) to zero at time $T$. This completes the proof of Theorem 1.1. \hfill$\Box$

\section{Proof of Theorem 1.2}
\setcounter{equation}{0}
In this section, we prove the uniqueness for our inverse source problem, i.e. Theorem 1.2.

\vspace{2mm}

{\noindent\bf Proof of Theorem 1.2.}\ 
Let $u=R_1 p$. By virtue of $u$ as a solution of equation \eqref{1.3}, we know that $p$ solves
\begin{equation}\label{5.1}
	\left\{
	\begin{aligned}
		&
		\begin{aligned}
		{\rm d}p-p_{xx}{\rm d}t-x^{2\gamma}p_{yy}{\rm d}t-\frac{\sigma}{x^2}p{\rm d}t =
		& \frac{2R_{1,x}}{R_1}p_x{\rm d}t + \frac{2x^{2\gamma}R_{1,y}}{R_1}p_y{\rm d}t \\
& \left( - \frac{R_{1,t}}{R_1} + \frac{R_{1,xx}}{R_1} + \frac{x^{2\gamma}R_{1,yy}}{R_1} \right)p{\rm d}t \\
		& +h(x,t){\rm d}t+\frac{R_2}{R_1}H(t){\rm d}B(t), \quad (x,y,t)\in Q_T,
		\end{aligned} \\
		& p(x,y,t)=0, \quad\quad\quad\qquad\qquad\qquad\qquad\qquad\qquad\qquad\qquad\qquad\quad\  (x,y,t)\in\Sigma_T,\\
        & p(x,y,0)=0, \quad\quad\quad\qquad\qquad\qquad\qquad\qquad\qquad\qquad\qquad\qquad\quad\  (x,y)\in I.
	\end{aligned}
	\right.
\end{equation}
Letting $w=p_y$, together with $u_y|_{\Sigma_T}=0$, $\mathbb P-a.s.$, we obtain
\begin{equation}\label{5.2}
	\left\{
	\begin{aligned}
		&
		\begin{aligned}
			{\rm d}w-w_{xx}{\rm d}t-x^{2\gamma}w_{yy}{\rm d}t-\frac{\sigma}{x^2}w{\rm d}t =&\frac{2R_{1,x}}{R_1}w_x{\rm d}t+\frac{2x^{2\gamma}R_{1,y}}{R_1}w_y{\rm d}t\\
& +\left(-\frac{R_{1,t}}{R_1}+\frac{R_{1,xx}}{R_1}+\frac{x^{2\gamma}R_{1,yy}}{R_1}\right)w{\rm d}t \\
			& +\left(\frac{2R_{1,x}}{R_1}\right)_y p_x{\rm d}t +\left(\frac{2x^{2\gamma}R_{1,y}}{R_1}\right)_yp_y{\rm d}t\\
&+\left(-\frac{R_{1,t}}{R_1}+\frac{R_{1,xx}}{R_1}+\frac{x^{2\gamma}R_{1,yy}}{R_1}\right)_y p{\rm d}t \\
          & +\left(\frac{R_2}{R_1}\right)_y H(t){\rm d}B(t),\quad\ \ \quad (x,y,t)\in Q_T,
		\end{aligned} \\
		& w(x,y,t)=0, \qquad\qquad\qquad\qquad\qquad\qquad\qquad\qquad\qquad\qquad\qquad (x,y,t)\in\Sigma_T,\\
        & w(x,y,0)=0, \qquad\qquad\qquad\qquad\qquad\qquad\qquad\qquad\qquad\qquad\qquad  (x,y)\in I.
	\end{aligned}
	\right.
\end{equation}
Applying Theorem 3.8 to $w$, we find that
\begin{align}\label{3-5.3}
&\mathbb{E}\int_{Q_T} s\lambda^2\Phi\Theta^2|w_x|^2{\rm d}x{\rm
d}y{\rm d}t+\mathbb{E}\int_{Q_T}s\lambda^2\Phi \Theta^2 x^{2\gamma}|w_y|^2{\rm d}x{\rm d}y{\rm d}t\nonumber\\
&+\mathbb{E}\int_{Q_T}s^3\lambda^4\Phi^3\Theta^2|w|^2{\rm d}x{\rm d}y{\rm d}t+\mathbb{E}\int_{Q_T}s\lambda\Phi \Theta^2\left|\left(\frac{R_2}{R_1}\right)_y\right|^2|H|^2{\rm d}x{\rm d}y{\rm d}t\nonumber\\
\leq &C\mathbb{E}\int_{Q_T}\Theta^2\left( |w_x|^2+x^{2\gamma}|w_y|^2+|w|^2+|p_x|^2+|p_y|^2+|p|^2\right){\rm d}x{\rm d}y{\rm d}t\nonumber\\
&+C\mathbb{E}\int_{Q_T}s\Phi\Theta^2\left|\nabla\left(\frac{R_2}{R_1}\right)_y\right|^2|H|^2{\rm d}x{\rm d}y{\rm d}t\nonumber\\
&+C\mathbb{E}\int_{I}s^2\lambda^2\Phi^2(T)\Theta^2(T)w^2(T){\rm d}x{\rm d}y+C\mathbb E\int_{\Gamma_T} s\lambda \Phi \Theta^2 |w_x|^2{\rm d}y{\rm d}t
\end{align}
for all $\lambda\geq\lambda_2$, $s\geq s_2$. By means of $w=p_y$ and $p(x,0,t)=0$ for $(x,t)\in I_x\times(0,T)$, we see that
\begin{align}
p(x,y,t)=\int_0^y w(x,\eta,t){\rm d}\eta.
\end{align}
Therefore, we obtain
\begin{align}\label{3-5.5}
&\mathbb E\int_{Q_T}\Theta^2\left(|p|^2+|p_x|^2+|p_y|^2\right){\rm d}x{\rm d}y{\rm d}t\leq C\mathbb E\int_{Q_T}\Theta^2\left(|w|^2+|w_x|^2\right){\rm d}x{\rm d}y{\rm d}t.
\end{align}
By (\ref{3-1.6}), we have
\begin{align}\label{3-5.6}
\mathbb{E}\int_{Q_T}s\Phi\Theta^2\left|\nabla\left(\frac{R_2}{R_1}\right)_y\right|^2|H|^2{\rm d}x{\rm d}y{\rm d}t\leq C\mathbb{E}\int_{Q_T}s\Phi\Theta^2\left|\left(\frac{R_2}{R_1}\right)_y\right|^2|H|^2{\rm d}x{\rm d}y{\rm d}t.
\end{align}
Thus, substituting (\ref{3-5.5}) and (\ref{3-5.6}) into (\ref{3-5.3}) and choosing $\lambda$ sufficiently large to absorb the first two terms on the right-hand side of (\ref{3-5.3}) by the terms on the left-hand side of (\ref{3-5.3}), we find that
\begin{align}\label{3-5.7}
&\mathbb{E}\int_{Q_T} s\lambda^2\Phi\Theta^2|w_x|^2{\rm d}x{\rm
d}y{\rm d}t+\mathbb{E}\int_{Q_T}s\lambda^2\Phi \Theta^2 x^{2\gamma}|w_y|^2{\rm d}x{\rm d}y{\rm d}t\nonumber\\
&+\mathbb{E}\int_{Q_T}s^3\lambda^4\Phi^3\Theta^2|w|^2{\rm d}x{\rm d}y{\rm d}t+\mathbb{E}\int_{Q_T}s\lambda\Phi\Theta^2\left|\left(\frac{R_2}{R_1}\right)_y\right|^2|H|^2{\rm d}x{\rm d}y{\rm d}t\nonumber\\
\leq&C\mathbb{E}\int_{I}s^2\lambda^2\Phi^2(T)\Theta^2(T)w^2(T){\rm d}x{\rm d}y+C\mathbb E\int_{\Gamma_T} s\lambda \Phi \Theta^2 |w_x|^2{\rm d}y{\rm d}t.
\end{align}

Since $u|_{\Gamma_T}=u_x|_{\Gamma_T}=0$, $\mathbb P$-a.s., we have $u_y|_{\Gamma_T}=u_{xy}|_{\Gamma_T}=0$ and further $w_x|_{\Gamma_T}=0$, $\mathbb P$-a.s.  Moreover $w(T)=0$ in $I$,   due to (\ref{3-11.8}). Then from (\ref{3-5.7}) we deduce
\begin{align}\label{5.6}
w=0 \quad {\rm in}\ Q_T,\quad \mathbb P-a.s.
\end{align}
which implies
\begin{align}\label{5.6}
u=0 \quad {\rm in}\ Q_T,\quad \mathbb P-a.s.
\end{align}
By (\ref{5.6}) and the equation (\ref{1.3}) of $u$,  we have
\begin{align}
\int_0^t h(x,\tau)R_1(x,y,\tau){\rm d}\tau +\int_0^t H(\tau )R_2(x,y,\tau){\rm d}B(\tau)=0,\quad t\in (0,T),
\end{align}
Together with (\ref{3-1.5}), we finally obtain (\ref{3-1.7}) and (\ref{3-1.8}). The proof of Theorem 1.2 is completed. \hfill$\Box$

\vspace{2mm}
\vskip0.5cm {\bf Acknowledgement.}\ This work is supported by NSFC
(No.11661004, No.11601240) \vskip 0.5cm

\newcounter{cankao}
\begin{list}
{[\arabic{cankao}]}{\usecounter{cankao}\itemsep=0cm} \centerline{\bf
References} \vspace*{0.5cm} \small

\vspace{2mm}
\item\label{Anh2016} C. T. Anh, V. M. Toi, Null controllability in large time of a parabolic equation involving
the Grushin operator with an inverse-square potential, Nonlinear Differential Equations and Applications 23 (2016) 1-26.

\item\label{Anh2013} C. T. Anh, V. M. Toi, Null controllability of a parabolic equation
involving the Grushin operator in some multi-dimensional domains,
Nonlinear Analysis: Theory, Methods and Applications 93 (2013) 181-196.

\item\label{Beauchard2014}  K. Beauchard, P. Cannarsa and R. Guglielmi,   Null
controllability of Grushin-type operators in dimension two.  Journal of the European Mathematical Society 16 (2014) 67-101.

\item\label{Beauchard2014IP} K. Beauchard, P. Cannarsa and M. Yamamoto, Inverse source problem and null
controllability for multidimensional parabolic operators of Grushin
type, Inverse Problems 30 (2014) 025006(26pp).

\item\label{Beauchard2015} K. Beauchard, L. Miller, M. Morancey, 2D Grushin-type equations:
minimal time and null controllable data, Journal of Differential Equations 259 (2015) 5813-5845.

\item\label{Barbu2003} V. Barbu, A. Rascanu, G. Tessitore, Carleman estimate and controllability of linear stochastic heat equations, Applied Mathematics and Optimization 47 (2003) 97-120.

\item\label{Buhgeim1981} A. Bukhgeim, M. V. Klibanov, Global uniqueness of a class of multidimentional inverse problems, Soviet Mathematics Doklady, 24 (1981), 244-247.

\item\label{Cannarsa2014} P. Cannarsa, R. Guglielmi, Null controllability in large time for the parabolic
Grushin operator with singular potential, Geometric Control Theory and
Sub-Riemannian Geometry, Springer International Publishing, (2014) 87-102.

\item\label{Cannarsa2005} P. Cannarsa, P. Martinze, J. Vancostenoble, Null controllability of degenerate heat equations, Advances in Differential Equations, 10 (2015) 153-190.

\item\label{Cannarsa2010} P. Cannarsa, J. Tort, M. Yamamoto, Determination of source terms in
a degenerate parabolic equation, Inverse Problems 26 (2010) 105003(26pp).


\item \label{Fragnelli2016} G. Fragnelli. Interior degenerate/singular parabolic equations in
nondivergence form: well-posedness and Carleman estimates, Journal of Differential Equations 260 (2016) 1314-1371.

\item\label{Fu2007} X. Fu, J. Yong, X. Zhang, Exact controllability for multidimensional semilinear hyperbolic equations, SIAM Journal on Control and Optimization 46 (2007) 1578-1614.

\item\label{Gao2014} P. Gao, Carleman estimate and unique continuation property for the linear stochastic Korteweg-de Vries equation, Bulletin of the Australian Mathematical Society 90 (02) (2014) 283-294.

\item\label{Gao2016} P. Gao,  A new global Carleman estimate for Cahn-Hilliard type
equation and its applications, Journal of Differential Equations 260 (2016) 427-444.

\item\label{Gao2015} P. Gao, M. Chen, Y. Li, Observability estimates and null controllability for forward and backward linear stochastic Kuramoto-Sivashinsky equations, SIAM Journal on Control and Optimization 53 (1) (2015) 475-500.

\item\label{Peng1991SAA} Y. Hu, S. Peng, Adapted solution of a backward semilinear stochastic evolution equations,
Stochastic Analysis and Applications 9 (1991) 445-459.
\item\label{Imanuvilov2003} O. Y. Imanuvilov, M. Yamamoto, Carleman inequalities for parabolic equations in Sobolev spaces of negative order and exact controllability for semilinear parabolic equations, Publications of the Research Institute for Mathematical Sciences 39 (2003) 227-274.

\item\label{Jiang2017} D. Jiang, Y. Liu, M. Yamamoto, Inverse source problem for the
hyperbolic equation with a time-dependent principal part, Journal of
Differential Equations 262 (2017) 653-681.

\item\label{Klibanov2004} M. V. Klibanov, A. Timonov, Carleman Estimates for Coefficient
Inverse Problems and Numerical Applications, VSP, Utrecht, 2004.

\item\label{Klibanov2013} M. V. Klibanov, Carleman estimates for global uniqueness,
stability and numerical methods for coefficient inverse problems, Journal of
Inverse Ill-Posed Problems 21 (2013) 477-560.

\item\label{Koenig} A. Koenig, Non null controllability of the Grushin equation in
2D, arXiv preprint arXiv:1701.06467 2017.

\item\label{Liu2014} X. Liu, Global Carleman estimate for stochastic parabolic equations and its application, ESAIM: Control, Optimisation and Calculus of Variations  20 (3) (2014) 823.

\item\label{LiuSIAM2019} X. Liu, Y. Yu, Carleman Estimates of Some Stochastic Degenerate Parabolic Equations and Application, SIAM Journal on Control and Optimization 57 (2019) 3527-3552.

\item\label{2012} Q. L\"{u}, Carleman estimate for stochastic parabolic equations and inverse stochastic parabolic problems, Inverse Problems 28 (4) (2012) 045008(18pp).

\item\label{2013} Q. L\"{u}, Observability estimate for stochastic Schr\"{o}dinger equations and its applications, SIAM Journal on Control and Optimization 51 (2013) 121-144.

\item \label{lu2013} Q. L\"{u}, Observability estimate and state observation problems for stochastic hyperbolic
equations, Inverse Problems 29 (2013) 095011.

\item \label{Lu2015CPAM} Q. L\"{u}, X. Zhang, Global uniqueness for an inverse stochastic hyperbolic problem with
three unknowns, Communications on Pure and Applied Mathematics 68 (2015) 948-63.

\item\label{Morancey2013} M. Morancey, About unique continuation for
a 2D Grushin equation with potential having an internal singularity,
arXiv preprint arXiv:1306.5616, 2013.


\item\label{Rousseau2012} J. L. Rousseau, G. Lebeau, On Carleman estimates for elliptic
and parabolic operators. Applications to unique continuation and control of parabolic equations, ESAIM: Control, Optimisation and Calculus of Variations 18 (2012) 712-747.

\item\label{Saut1987} J. C. Saut, B. Scheurer, Unique continuation for some
evolution equations, Journal of Differential Equations 66 (1987) 118-139.

\item\label{Zhang2009} S. Tang, X. Zhang, Null controllability for forward and backward stochastic parabolic equations, SIAM Journal on Control and Optimization 48 (2009) 2191-2216.

\item\label{Wang2014} C. Wang, R. Du, Carleman estimates and null controllability for a class of degenerate parabolic equations with convection terms, SIAM Journal on Control and Optimization 52 (2014) 1457-1480.

\item\label{Wu2017} B. Wu, J. Yu, H\"{o}lder stability of an inverse problem for a
strongly coupled reaction-diffusion system, IMA Journal of Applied Mathematics 82 (2017) 424-444.

\item\label{Wu2019JIIP} B. Wu, Y. Gao, Z. Wang and Q. Chen,  Unique continuation for a reaction-diffusion
system with cross diffusion, Journal of Inverse and Ill-posed Problems 27 (2019) 511-525.

\item\label{Yan2018JMAA} Y. Yan, Carleman estimates for stochastic parabolic equations with Neumann boundary conditions and applications, Journal of  Mathematical Analysis and Applications 457 (2018) 248-272.

\item\label{Yamamoto2009} M. Yamamoto, Carleman estimates for parabolic equations and
applications, Inverse Problems 25 (2009) 123013(75pp).

\item\label{Yuan2015} G. Yuan, Determination of two kinds of sources simultaneously for a stochastic wave
equation, Inverse Problems 31 (2015) 085003(13pp).

\item\label{Zhang2008} X. Zhang, Carleman and observability estimates for stochastic wave equations, SIAM Journal on Mathematical Analysis 40 (2008) 851-868.

\item\label{Zhou1992} X. Zhou, A duality analysis on stochastic partial differential equations, Journal of Functional Analysis 103 (1992) 275-193.







\end{list}

\end{document}